\documentstyle{amsppt}

\magnification=1100 \NoBlackBoxes\nologo %\hsize=11.5cm

\def\half{\frac{1}{2}}

\def\inv{^{-1}}
\def\?{{\bf{??}}}

\def\Hilb{\text{Hilb}}

\def\Proj{\text{Proj}}

\def\A{\Bbb A}

\def\C{\Bbb C}

\def\P{\Bbb P}
\def\R{\Bbb R}

\def\ord{\text{\rm ord}}

\def\Spec{\text{\rm Spec} }

\def\Q{\Bbb Q}
\def\[{\big[}
\def\]{\big]}

\def\O{\Cal O}

\def\y{\bar{y}}

\def\Sym{\text{Sym}}

\def\rk{\text{rk}}

\def\m{\frak m}

\def\1/2{\frac{1}{2}}

\def\I{\Cal I}

\def\2{{[2]}}
\def\l{\ell}
\def\Proj{\text{Proj}}
\def\sbr #1.{^{[#1]}}
\def\scl #1.{^{\lceil#1\rceil}}
\def\spr #1.{^{(#1)}}

\def\nl{\newline}
%\input myheader.tex
%\nologo
\topmatter
\title Geometry on nodal curves II:
\centerline{ cycle map and intersection calculus}
\endtitle
\author
Ziv Ran
%\thanks{\raggedright{
%Partially supported by NSA Grant MDA904-02-1-0094} }
\endauthor
\date 2004.09.20\enddate

\address University of California, Riverside CA 92521\endaddress
\email ziv\@math.ucr.edu\endemail 
\endkeywords
%\urladdr http://www.math.ucr.edu/~ziv/papers/
%semicurv.pdf\endurladddr
\rightheadtext { Cycle Map and Intersection Calculus}
\leftheadtext{Z. Ran} \abstract We study the relative Hilbert scheme
of
 a family of nodal (or smooth)
curves via its (birational) {\it{ cycle map}}, going to the relative
symmetric product.  We show the cycle map is the blowing up of the
discriminant locus, which consists of cycles with multiple points.
We derive an intersection calculus for Chern classes of tautological
bundles on the relative Hilbert scheme, which has applications to
enumerative geometry.

\endabstract
 \thanks \raggedright {Updates and corrections
available at $\underline{math.ucr.edu/\tilde{\ }
ziv/papers/cyclemap.pdf}$}\linebreak Research Partially supported
by NSA Grant MDA904-02-1-0094; reproduction and distribution of
reprints by US government permitted.
\endthanks

\endtopmatter\document
Consider a family of curves given by a flat projective morphism
$$\pi:X\to B$$
 over an irreducible (and usually
projective) base, with fibres $$X_b=\pi\inv(b), b\in B$$ which are
irreducible  nonsingular for the generic $b$ and   at worst nodal
for every $b$. Many questions in the classical projective and
enumerative geometry of this family can be naturally  phrased, and
in a formal sense solved (see for instance \cite{R}), in the context
of the {\it{relative Hilbert scheme}}
$$X\sbr m._B=\Hilb_m(X/B),$$ which parametrizes length-$m$
subschemes of $X$ contained in fibres of $\pi$, and the natural {\it
{tautological vector bundles}} that live on $X\sbr m._B$. Typically,
the questions include ones involving relative multiple points and
multisecants in the family, and the formal solutions involve Chern
numbers of the tautological bundles. Thus, turning these formal
solutions into meaningful ones requires computing the Chern numbers
in question.\par This paper is a contribution to the study, both
qualitative and enumerative, of the relative Hilbert scheme of a
family of modal curves as above. We provide the following:\nl -- a
structure theorem for the cycle (or 'Hilb-to-Chow') map
$$\frak c_m:X\sbr m._B\to X\spr m._B,$$ where $X\spr m._B$ is the
relative symmetric product, showing that $\frak c_m$ is equivalent
to the blowing up of the {\it{discriminant locus}} $$D^m\subset
X\spr m._B,$$ which parametrizes nonreduced cycles;\nl -- when $X$
is a smooth surface, an intersection calculus for certain
'tautological classes' allowing computation of the Chern numbers of
the tautological bundles on $X\sbr m._B$.\par To be precise, this
calculus, which is based on the structure theorem, actually takes
place on the (full) flag-Hilbert scheme $W^m(X/B)$, parametrizing
length-$m$ flags of subschemes of fibres of $X/B$, whose basic
theory was developed in \cite{R}. Nonethless,  the Chern numbers
computed in this calculus are the same, up to an evident factor, as
those on $X\sbr m._B$. Using the calculus, it is possible to compute
explicitly the expressions given in \cite{R} for various
multiple-point and multisecant cycles. The advantage of using
$W^m(X/B)$ over $X\sbr m._B$ is that the tautological classes are
expressed as polynomials in  {\it divisor} classes $\Gamma\scl i.,
i=2,...,m$, corresponding to certain diagonal loci, together with
the classes coming from $X$ itself. This allows us to work in the
ring $T^m$ generated by these classes, a ring that we call the
{\it{tautological ring}} on $W^m(X/B)$. Working in $T^m$, one is
effectively working with divisor classes-- in fact, $T^m$ contains
explicit expressions for the {\it{Chern roots}} of the tautological
bundles, which are convenient in computations. Thus, passage to
$W^m(X/B)$ and its tautological ring may be viewed as a version of
the familiar 'splitting principle'.\par What our calculus does is,
essentially, to compute the operator of multiplication by
$\Gamma\scl m.$ on $T^m.$
 To be precise, our
method effectively yields a set of additive generators of $T^m$,
together with rules for expressing the product of a generator with
$\Gamma\scl m.$ as linear combination of generators. Given the
inductive structure in $m$ of the $T^m$, this completely determines
the ring structure on $T^m$, albeit with an apparent ambiguity if
(and only if) our generators are linearly dependent. It seems
reasonable to conjecture that our generators are in fact linearly
independent, but we do not prove this. In any event, our calculus is
certainly sufficient to compute the top-degree products, which are
those with enumerative significance.
\par Note that if $X$ is a
smooth surface, there is a natural closed embedding
$$j_\pi\sbr m.:X\sbr m._B\subset X\sbr m.$$ of the relative Hilbert scheme in
the full Hilbert scheme of $X$, which is a smooth projective
$2m$-fold. There is a large literature on Hilbert schemes of smooth
surfaces and their cohomology and intersection theory, due to
Ellingsrud-Str{\o}mme, G\"ottsche, Nakajima, Lehn and others, see
\cite{EG, L, LS, N} and references therein. In particular, Lehn
\cite{L} gives a formula for the Chern classes of the tautological
bundles on the full Hilbert scheme $X\sbr m.$, from which one can
derive a formula for the analogous classes on $X\sbr m._B$ if $X$ is
a smooth surface, but this does not, to our knowledge, yield Chern
numbers (besides the top one) on $X\sbr m.,$ much less $X\sbr m._B$
(the two sets of numbers are of course different). Going from Chern
{\it classes} to Chern {\it numbers} it a matter of working out the
top-degree multiplicative structure, i.e. the intersection calculus.
 When $X$ is a surface with trivial
canonical bundle, Lehn and Sorger \cite{LS} have given a rather
involved description of the mutiplicative structure on the
cohomology of $X\sbr m.$ in all degrees, not just the top one. While
products on $X\sbr m.$ and $X\sbr m._B$ are compatible $j_\pi\sbr
m.$, it's not clear how to compute intersection products, especially
intersection {\it numbers} on $X\sbr m._B$ from products on $X\sbr
m.$, even in case $X$ has trivial canonical bundle. Indeed some of
our additive generators directly involve the {\it{fibre nodes}} of
the family $X/B$ and do not appear to come from classes on $X\sbr
m.$. In any event, the computing the relative cohomology of the pair
$(X\sbr m._B, X\sbr m.)$ is an interesting problem that at the
moment seems out of reach.
\par

\subheading{Acknowledgements} I thank Mirel Caibar  and Rahul
Pandharipande for valuable comments. A preliminary version of some
of these results was presented at conferences in Siena, Italy and
Hsinchu, Taiwan, in June 2004, and I thank the organizers of these
conferences for this opportunity.

\heading 0. Preliminaries\endheading We define a combinatorial
function that will be important in computations to follow. Denote
by $Q$ the closed 1st quadrant in the real $(x,y)$ plane,
considered as an additive cone. We will consider unbounded
$Q$-invariant closed subsets $R\subset Q$ with the property that
the boundary of $R$ relative to $Q$ consists of a finite number of
finite horizontal and vertical segments with integral endpoints
(the boundary of $R$ in $\R^2$ will then consist of this plus two
semi-infinite intervals, one on each axis).
% , plus 2
%semi-infinite intervals, one each on the $x$ and $y$ axes.
We call such $R$ a {\it{special infinite polygon}}. The closure of
the complement $$S=R^c:=\overline{Q\setminus R}\subset Q$$ has
finite (integer) area and will be called a {\it{special finite
polygon}}; in fact the area of $S$ coincides with the number of
integral points in $S$ that are $Q$-interior, i.e. not in $R$;
these are precisely the integer points $(a,b)$ such that
$[a,a+1]\times [b,b+1]\subset S$. Fixing a natural number $m$, the
basic special finite polygon associated to $m$ is
$$S_m=\bigcup_{i=1}^m
[0,\binom{m-i+1}{2}]\times[0,\binom{i+1}{2}].$$ It has area
$$\alpha_m=\sum\limits_{i=1}^{m-1} i\binom{m+1-i}{2}
=3\binom{m}{4}+3\binom{m}{3}+m-1$$ and associated special infinite
polygon denoted $R_m.$ Now for each integer $j=1,...,m-1$ we
define a special infinite polygon $R_{m,j}$ as follows. Set
$$P_j=(-j, m+1-j),$$
$$R_{m,j}=\left( R_m\cup(R_m+P_j)\cup[0,\infty)\times[j,\infty)
\right )\cap Q$$ (where $R_m+P_j$ denotes the translate of $R_m$ by
$P_j$ in $\R^2$). Then let $S_{m,j}=R_{m,j}^c,$
$$\beta_{m,j}={\text {area}}(S_{m,j}),$$
$$\beta_m=\sum_{j=1}^{m-1}\beta_{m,j}.$$ It is easy to see that
$$\beta_{m,1}=\binom{m}{2}, \beta_{m,j}=\beta_{m,m-j}$$ but
otherwise we don't know a closed-form formula for these numbers in
general. A few small values are
$$\beta_{2,1}=\beta_2=1$$
$$\vec{\beta}_{3}=(3,3), \beta_3=6$$ $$\vec{\beta}_{4}=(6,8,6),
\beta_4=20$$ $$\vec{\beta}_{5}=(10, 15,15,10), \beta_5=50$$
$$\vec{\beta}_{6}=(15,24,27,24,15), \beta_6=105.$$ For an interpretation
for these numbers see \S 1.6 below.
\comment
is as follows. Let
$\theta_m$ be the quotient of the $\C$-algebra $ \C[x,y]$ by the
ideal
$$J_m=(x^\binom{m}{2},...,x^\binom{m-i+1}{2}y^\binom{i}{2},...,
y^\binom{m}{2}),$$ so $\Spec (\theta_m)$ is a subscheme of
$\A^2_\C$ supported at the origin. Then  $\theta_m$ has a $\C$-
basis of the form $\{x^\mu
y^\nu:(\mu,\nu)\in{\text{interior}}_Q(S_m)\}$ and in particular
$\dim_\C(\theta_m)=\alpha_m.$ Let $\theta_{m,j}$ be the quotient
of $\theta_m$ by the principal ideal
$$J_{m,j}=(y^j+\eta x^{m-j})$$ where $\eta\in\C$ is any nonzero number.
Then it
is easy to see that $$\beta_{m,j}=\dim_\C(\theta_{m,j}).\tag ?$$
We will use the $\beta_{m,j}$ through this equality.
\endcomment

 \heading \bf{1. The cycle map as
blowup}\endheading\subheading{1.1 Set-up} Let
$$\pi : X\to B\tag 1.1.1$$ be a family of nodal (or smooth) curves
with $X, B$ smooth. Let $X^m_B,  X\spr m._B$, respectively, denote
the $m$th Cartesian and symmetric fibre products of $X$ relative
to $B$. Thus, there is a natural map
$$\omega_m:X^m_B\to  X\spr m._B\tag 1.1.2$$ which realizes its target
as the quotient of its source under the permutation action of the
symmetric group $\frak S_n.$ Let $$\Hilb_m(X/B)=X\sbr m._B$$
denote the relative Hilbert scheme paramerizing length-$m$
subschemes of fibres of $\pi$, and
$$\frak c= \frak c_m :X\sbr m._B\to  X\spr m._B\tag 1.1.3$$ the natural
{\it{cycle map}} (cf.\cite{A}). Let $D^m\subset X\spr m._B$ denote
the discriminant locus or 'big diagonal', consisting of cycles
supported on $<m$ points (endowed with the reduced scheme
structure). Clearly, $D^m$ is a prime Weil divisor on $ X\spr
m._B$, birational to $X\times_B\Sym^{m-2}(X/B)$, though it is less
clear what the defining equations of $D^m$ on $ X\spr m._B$ are
near singular points. The purpose of this section is to prove
\proclaim{Theorem 1} The cycle map $$\frak c_m:X\sbr m._B\to
 X\spr m._B$$ is the blow-up of
$D^m\subset X\spr m._B$.\endproclaim \subheading{1.2 Preliminary
reductions} To begin with, we reduce the Theorem to a local
statement over a neighborhood of a 1-point cycle $mp\in X\spr m._B$
where $p\in X$ is a node of $\pi\inv(\pi(p))$. Set
$$\Gamma\spr m.=\frak c_m\inv(D^m)\subset X\sbr m._B.\tag 1.2.1 $$
It was shown in
\cite{R}, and will be reviewed below, that $\frak c_m$ is a small
birational map (with fibres of dimension $\leq 1$), and that
$X\sbr m._B$ is smooth. Consequently $\Gamma\spr m.$ is an
integral, automatically Cartier, divisor, and therefore $\frak c$
factors through a map $\frak c'$ to the blow-up $B_{D^m}( X\spr
m._B)$, and it suffices to show that $\frak c'$ is an isomorphism,
which can be checked locally.
\par Next, let $X^o\subseteq X$ denote the open subset
consisting of regular points of $\pi$, i.e. points $x\in X$ where
$\pi$ is smooth (submersive) or equivalently, such that $x$ is a
smooth point of $\pi\inv(\pi(x))$. Note that the open subset
$\Sym^m(X^o/B)\subseteq X\spr m._B$ is smooth and
$$\frak c_m:\frak c_m\inv(\Sym^m(X^o/B))\to\Sym^m(X^o/B)$$
is an isomorphism. Therefore it will suffice to show $\frak c_m$
is equivalent to the blowing-up of $D^m$ locally near any cycle
$Z\in X\spr m._B$ whose support meets the locus $X^\sigma\subset
X$ of singular points of $\pi$ (i.e. singular points of fibres).
Writing
$$Z=\sum \limits_{i=1}^km_ip_i$$ with $m_i>0, p_i$ distinct, we
have a cartesian diagram $$\matrix \prod\limits_{i=1}^k {_B}\
X\sbr m_i._B&\overset \prod \frak c_{m_i}\to \longrightarrow&
\prod\limits_{i=1}^k {_B}\ X\spr m_i._B\\e_1\uparrow&&\ \ \
\uparrow d_1\\
\ \ \ H &\to&S\\
e\ \downarrow &&\ \ \downarrow d\\
X\sbr m._B&\overset \frak c_m\to\longrightarrow & X\spr m._B
\endmatrix\tag 1.2.2
$$ Where $H$ is the natural inclusion correspondence on Hilbert
schemes: $$H=\{ (\zeta_1,...,\zeta_k, \zeta)\in
\prod\limits_{i=1}^k {_B}\ X\sbr m_i._B\times X\sbr
m._B:\zeta_i\subseteq\zeta, i=1,...,k\},$$ and similarly for $S$.

  Note that the
right vertical arrows $d, d_1$ are isomorphisms between some
neighborhoods $U$ of $Z$ and $U'$ of $(m_1p_1,...,m_kp_k)$ and the
left vertical arrows $e, e_1$ are isomorphisms between $\frak
c_m\inv(U)$ and $(\prod \frak c_{m_i})\inv(U')$.  Now by
definition, the blow-up of $ X\spr m._B$ in $D^m$ is the Proj of
the graded algebra
$$A(\I_{D^m})=\bigoplus\limits_{n=0}^\infty\I_{D^m}^{\ n}.$$ Note that
$$d\inv(D^m)=\sum p_i\inv(D^{m_i})$$ and moreover,
$$d^*(\I_{D^m})=\bigotimes{_B}\ p_i^*(\I_{D^{m_i}})$$ where we use
$p_i$ generically to denote an $i$th coordinate projection.
Therefore,
$$A(\I_{D^m})\simeq\bigotimes{_B}\ p_i^*A(\I_{D^{m_i}})$$ as graded
algebras,  compatibly with the isomorphism
$$\O_{\prod\limits_{i=1}^k {_B}\ \Sym^{m_i}(X/B)}\simeq\bigotimes
\limits_{i=1}^k{_B}\  \O_{ \Sym^{m_i}(X/B)}.$$ Now it is a general
fact that Proj is compatible with tensor product of graded
algebras, in the sense that $$\Proj(\bigotimes {_B}\ A_i)\simeq
\prod {_B}\ \Proj(A_i).$$ Consequently (1.2.2) induces another
cartesian diagram with unramified vertical arrows
$$\matrix \prod\limits_{i=1}^k {_B}\
X\sbr m_i._B&\overset \prod c'_{m_i}\to \longrightarrow&
\prod\limits_{i=1}^k {_B}\ B_{D^{m_i}}X\spr m_i._B\\
\downarrow&&\downarrow \\
X\sbr m._B&\overset c'_m\to\longrightarrow &B_{D^m} X\spr m._B.
\endmatrix\tag 1.2.3
$$
To prove $c'_m$ is an isomorphism, it suffices to prove that so is
$c'_{m_i}$ for each $i$. \comment
 an analytic
neighborhood of $Z$ in $ X\spr m._B$ has the form
$\prod\limits_{i=1}^k {_B}\Sym^{m_i}_B(U_i)$, with $U_i$ a
neighborhood of $p_i$, and an analogous cartesian product
decomposition holds locally also for $\Hilb_m$ and for the map
$c$.
\endcomment
The upshot of this is that it suffices to prove $c=\frak c_m$ is
equivalent to the blow-up of $ X\spr m._B$ in $D^m$, locally over
a neighborhood of a cycle of the form $mp$ where $p\in X$ is a
singular point of $\pi.$\par \subheading{1.3 A local model} Fixing
such a point $p$, we have coordinates on an affine neighborhood
$U$ of $p$ in $X$ so that $\pi$ is given on $U$ by
$$t=xy.$$ Then the relative cartesian product $X^m_B$, as
subscheme of $X^m\times B$, is given by
$$x_1y_1=...=x_my_m=t.\tag 1.3.1$$ Let $\sigma_i^x, \sigma_i^y,
i=0,...,m$ denote the elementary symmetric functions in
$x_1,...,x_m$ and in $y_1,...,y_m$, respectively, where we set
$\sigma_0=1$. Put together with the projection to $B$, they yield
a map
$$\sigma:\Sym^m(U/B)\to \A^{2m}_B= \A^{2m}\times B\tag 1.3.2$$
 $$\sigma=((-1)^m\sigma_m^x,...,-\sigma_1^x,
(-1)^m\sigma_m^y,...,-\sigma_1^y, \pi^{(m)})$$ where $\pi^{(m)}:
X\spr m._B\to B$ is the structure map.
  \proclaim{Lemma 2} $\sigma$ is an
embedding locally near $mp$.\endproclaim\demo{proof} It suffices
to prove this formally, i.e. to show that $\sigma_i^x, \sigma_j^y,
i,j=1,...,m$ generate topologically the completion $\hat{\m}$ of
the maximal ideal of $mp$ in $ X\spr m._B.$ To this end it
suffices to show that any $\frak S_m$-invariant polynomial in the
$x_i, y_j$ is a polynomial in the $\sigma_i^x, \sigma_j^y$ and
$t$. Let us denote by $R$ the averaging or symmetrization operator
with respect to the permutation action of $\frak S_m$, i.e.
$$R(f)=\frac{1}{m!}\sum\limits_{g\in\frak S_m}g^*(f).$$
 Then it suffices to show that the elements
$R(x^Iy^J)$, where $x^I$ (resp. $y^J$) range over all monomials in
$x_1,...,x_m$ (resp. $y_1,...,y_m$) are polynomials in the
$\sigma_i^x, \sigma_j^y$ and $t$. Now the relation (1.3.1)
defining $X^m_B$ easily implies that
$$R(x^Iy^J)-R(x^I)R(y^J)=tF$$ where $F$ is an $\frak
S_m$-invariant polynomial in the $x_i, y_j$ of bidegree
$(|I|-1,|J|-1)$, hence a linear combination of elements of the
form $R(x^{I'}y^{J'}), |I'|=|I|-1, |J'|=|J|-1$. By induction, $F$
is a polynomial in the $\sigma_i^x, \sigma_j^y$ and clearly so is
$R(x^I)R(y^J).$ Hence so is $R(x^Iy^J)$ and we are
done.\qed\enddemo

Now let $C_1,...,C_{m-1}$ be copies of $\P^1$, with homogenous
coordinates $u_i,v_i$ on the $i$-th copy. Let $\tilde{C}\subset
C_1\times...\times C_{m-1}\times B$ be the subscheme defined by
$$v_1u_2=tu_1v_2,..., v_{m-2}u_{m-1}=tu_{m-2}v_{m-1}.\tag 1.3.3$$ Thus
$\tilde{C}$ is a reduced complete intersection of divisors of type
$(1,1,0,...,0), (0,1,1,0,...,0)$ ,..., $(0,...,0,1,1)$ and it is
easy to check that the fibre of $\tilde{C}$ over $0\in B$ is
$$\tilde{C}_0=\bigcup\limits_{i=1}^m\tilde{C}_i,$$
where $$\tilde{C}_i= [1,0]\times...\times[1,0]\times
C_i\times[0,1]\times...\times[0,1]$$ and that in a neighborhood of
$\tilde{C}_0$, $\tilde{C}$ is smooth and $\tilde{C}_0$ is its
unique singular fibre over $B.$ We may embed $\tilde{C}$ in
$\P^{m-1}\times B,$ relatively over $B$ using the mutihomogenous
monomials $$Z_i=u_1\cdots u_{i-1}v_{i}\cdots v_{m-1}, i=1,...,m.$$
These satisfy the relations $$Z_iZ_j=tZ_{i+1}Z_{j-1}, i<j-1\tag
1.3.4$$ so they embed $\tilde{C}$ as a family of rational normal
curves $\tilde{C}_t\subset\P^{m-1}, t\neq 0$ specializing to
$\tilde{C}_0$, which is embedded as a nondegenerate, connected
$(m-1)$-chain
 of lines.\par
 Next consider an affine
space $\A^{2m}$ with coordinates $a_0,...,a_{m-1},
d_0,...,d_{m-1}$ and let $\tilde{H}\subset\tilde{C}\times\A^{2m}$
be the subscheme defined by
$$a_0u_1=tv_1, d_0v_{m-1}=tu_{m-1}$$
$$a_1u_1=d_{m-1}v_{1},...,a_{m-1}u_{m-1}=d_1v_{m-1}.\tag 1.3.5$$
Set $L_i=p_{C_i}^*\O(1).$ Then consider the subscheme of
$Y=\tilde{H}\times_{B}U$ defined by the equations
$$F_0:=x^m+a_{m-1}x^{m-1}+...+a_1x+a_0\in
\Gamma(Y,\O_Y)$$
$$F_1:=u_1x^{m-1}+u_1a_{m-1}x^{m-2}+...+u_1a_2x+u_1a_1+v_1y
\in\Gamma(Y,L_1)$$ ...
$$F_i:=u_ix^{m-i}+u_ia_{m-1}x^{m-i-1}+...+u_ia_{i+1}x+u_ia_i+
v_id_{m-i+1}y+...+v_id_{m-1}y^{i-1}+ v_iy^i$$
$$
\in\Gamma(Y,L_i)\tag 1.3.6$$ ...
$$F_m:=d_0+d_1y_1+...+d_{m-1}y^{m-1}+y^m\in
\Gamma(Y,\O_Y).$$ The following statement summarizes results from
\cite{R1} \proclaim{Theorem 3} (i) $\tilde{H}$ is smooth and
irreducible.\par (ii) The ideal sheaf $\I$ generated by
$F_0,...,F_m$ defines a subscheme of $\tilde{H}\times_BX$ that is
flat of length $m$ over $\tilde{H}$
\par (iii)The classifying map $$\Phi=\Phi_\I:
\tilde{H}\to\Hilb_m(U/B)$$
is an isomorphism.
\endproclaim\demo{proof} The smoothness of $\tilde{H}$ is clear
from the defining equations equations and also follows from
smoothness of $\Hilb_m(U/B)$ once (ii) and (iii) are proven. To
that end consider the point $q_i, i=1,...,m,$ on the special fibre
of $\tilde{H}$ over $\A^{2m}_B$ with coordinates $$v_j=0,\ \forall
j< i; u_j=0,\ \forall j\geq i.$$ Then $q_i$ has an affine
neighborhood $U_i$ in $\tilde{H}$ defined by $$U_i=\{ u_j=1, \
\forall j< i;\ v_j=1, \ \forall j\geq i\},\tag 1.3.7$$ and these
$U_i, i=1,..., m$ cover a neighborhood of the special fibre of
$\tilde{H}.$ Now the generators $F_i$ admit the following
relations:
$$u_{i-1}F_j=u_jx^{i-1-j}F_{i-1},\ 0\leq j<i-1;\ v_iF_j=v_jy^{j-i}F_i,\
m\geq j>i$$ where we set $u_i=v_i=1$ for $i=0,m.$ Hence $\I$ is
generated there by $F_{i-1}, F_i$ and assertions (ii), (iii)
follow directly from Theorems 1,2 and 3 of \cite{R1}.\qed

\enddemo
\remark{Remark 3.1} For future reference, we note that over $U_i,$ a
co-basis for the universal ideal $\I$ (i.e. a basis for $\O/\I$) is
given by $1,...,x^{m-i}, y,...,y^{i-1}$. In view of the definition
of the $F_i$ (1.3.6), this is immediate from the fact just noted
that, over $U_i,$ the ideal $\I$ is generated by $F_{i-1}, F_i$,
plus the fact that on $U_i$ we have
$u_{i-1}=v_{i}=1.$\qed\endremark\remark{Remark 1.3.2} For integers
$\alpha, \beta\leq m$, consider the locus $X\spr m._B(\alpha,
\beta)$ of cycles containing $\alpha p+y'+y"$ where $p$ is a node
and  $y', y"$ are general cycles of degree $\beta$ (resp.
$m-\alpha-\beta$) on the two (smooth) components of the special
fibre. Then it is easy to see that the general fibre of  $\frak c_m$
over $X\spr m._B(\alpha)$ coincides with
$\bigcup\limits_{i=m-\beta-\alpha+1}^{m-\beta-1} C^m_i$, which may
be naturally identified with $C^{\alpha}=\bigcup\limits
_{j=1}^{\alpha-1}C^{\alpha}_j.$

\endremark

\subheading{1.4 Reverse engineering}
 In light of Theorem 3, we
identify a neighborhood $H_m$ of the special fibre in $\tilde{H}$
with a neighborhood of the punctual Hilbert scheme (i.e. $\frak
c_m\inv(mp)$) in $X\sbr m._B$, and note that the projection
$H_m\to \A^{2m}\times B$ coincides generically, hence everywhere,
with $\sigma\circ \frak c_m$. Hence $H_m$ may be viewed as the
subscheme of $\Sym^m(U/B)\times_B\tilde{C}$ defined by the
equations
$$\sigma_m^xu_1=tv_1,$$$$ \sigma_{m-1}^xu_1=\sigma_{1}^yv_{1},...,
%d_0v_{m-1}=tu_{m-1}$$
\sigma^x_{1}u_{m-1}=\sigma^y_{m-1}v_{m-1},\tag 1.4.1$$
$$tu_{m-1}=\sigma^y_mv_{m-1}$$
Alternatively, $H_m$ may be defined as the subscheme of
$\Sym^m(U/B)\times \P^{m-1}\times B$ defined by the relations
(1.3.3), which define $\tilde{C}$, together with
$$\sigma_{m-j}^yZ_i=t^{m-j-i}\sigma_j^xZ_{i+1},\ \
 i=1,...,m-1,j=0,...,m-1; \tag 1.4.2$$
$$\sigma_{m-j}^xZ_i=t^{m-j-i}\sigma_j^yZ_{i-1},\ \
i=2,...,m, j=0,...,m-1. \tag 1.4.3$$ Our task now is effectively
to 'reverse-engineer' an ideal in the $\sigma$'s whose syzigies
are given by (1.4.2-1.4.3). To this end, it is convenient to
introduce order in the coordinates. Thus let
$OH_m=H_m\times_{\Sym^m(U/B)}U^m_B$, so we have a cartesian
diagram $$\matrix OH_m&\overset \varpi_m\to\longrightarrow&
H_m\\ o\frak c_m \downarrow &&\downarrow \frak c_m\\
X^m_B&\overset \omega_m\to\longrightarrow& X\spr
m._B\endmatrix%\tag 1.4.4
$$ and its global analogue
$$\matrix X^{\lceil m\rceil}_B&\overset \varpi_m\to\longrightarrow&
X\sbr m._B\\ o\frak c_m \downarrow &&\downarrow \frak c_m\\
X^m_B&\overset \omega_m\to\longrightarrow& X\spr m._B\endmatrix\tag
1.4.4
$$

Note that $ X\spr m._B$ is normal, Cohen-Macaulay and
$Q$-Gorenstein: this follows from the fact that it is a quotient by
$\frak S_m$ of $X^m_B$, which is a locally complete intersection
with singular locus of codimension $\geq 2$ (in fact, $>2$, since
$X$ is smooth). Alternatively, normality of $ X\spr m._B$ follows
from the fact that $H_m$ is smooth and the fibres of \nl$\frak
c_m:H_m\to
 X\spr m._B$ are connected (being products of connected chains of
rational curves). Note that $\omega_m$ is simply ramified
generically over $D^m$ and we have $$\omega_m^*(D^m)=2OD^m $$ where
$$OD^m=\sum\limits_{i<j}D^m_{i,j}$$ where
$D^m_{i,j}=p_{i,j}\inv(OD^2)$ is the locus of points whose $i$th and
$j$th components coincide.
 To prove $\frak c_m$ is equivalent to the
blowing-up of $D^m$ it will suffice to prove that $o\frak c_m$ is
equivalent to the blowing-up of $2OD^m=\omega_m^* (D^m)$ which in
turn is equivalent to the blowing-up of $OD^m.$ \comment Note that
$OD^m$ is reduced but, unlike $D^m$, splits up as
$$OD^m=\sum\limits_{i<j}D^m_{i,j}$$ where
$D^m_{i,j}=p_{i,j}\inv(OD^2)$ is the locus of points whose $i$th
and $j$th components coincide.
\endcomment
The advantage of working with $OD^m$ rather than its unordered
analogue is that at least some of its equations are easy to write
down: let $$v^m_x=\prod_{1\leq i<j\leq m}(x_i-x_j),$$ and likewise
for $v^m_y.$ As is well known, $v^m_x$ is the determinant of the Van
der Monde matrix
$$V^m_x=\left [\matrix 1&\ldots&1\\x_1&\ldots&x_m\\\vdots&&\vdots\\
x_1^{m-1}&\ldots&x_m^{m-1}\endmatrix \right ].$$
 Also
set
$$\tilde{U}_i=\varpi_m\inv(U_i),$$ where $U_i$ is as in (1.3.7),
being a neighborhood of $q_i$ on $H_m$. Then in $U_1$, the
universal ideal $\I$ is defined by
$$F_0, \ \ F_1=y+(\text{function of }x)$$ and consequently the
length-$m$ scheme corresponding to $\I$ maps isomorphically to its
projection to the $x$-axis. Therefore over
$\tilde{U}_1=\varpi_m\inv(U_0),$ where $F_0$ splits as $\prod
(x-x_i),$ the equation of $OD^m$ is simply given by $$G_1=v^m_x.$$
Similarly, the equation of $OD^m$ in $\tilde{U}_m$ is given by
$$G_m=v^m_y.$$ New let $$\Xi:OH_m\to\P^{m-1}$$ be the morphism
corresponding to $[Z_1,...,Z_m]$, and set $L=\Xi^*(\O(1)).$ Note
that $\tilde{U}_i$ coincides with the open set where $Z_i\neq 0$,
so $Z_i$ generates $L$ over $\tilde{U}_i.$ Let
$$O\Gamma\spr m.=o\frak c_m\inv(OD^m).$$
This is a $1/2-$Cartier divisor because $2O\Gamma\spr
m.=\varpi_m\inv(\Gamma\spr m.)$ and $\Gamma\spr m.$ is Cartier,
$H_m$ being smooth. In any case, the ideal $\O(-O\Gamma\spr m.)$ is
a divisorial sheaf (reflexive of rank 1).
 Our aim is to construct an isomorphism
$$\gamma:\O(-O\Gamma\spr m.)\to L.\tag 1.4.5$$
Since $L=\Xi^*(\O(1))$ and
$OH_m$ is a subscheme of $X^m_B\times\P^{m-1},$ this isomorphism
would clearly imply Theorem 1. To construct $\gamma,$ it suffices
to specify it on each $\tilde{U}_i.$ \subheading{1.5 Mixed Van der
Mondes and conclusion of proof}  A clue as to how this might be
done comes from the relations (1.4.2-1.4.3). Thus, set
$$G_i=\frac{(\sigma_m^y)^{i-1}}{t^{(i-1)(m-i/2)}}v^m_x=
\frac{(\sigma_m^y)^{i-1}}{t^{(i-1)(m-i/2)}}G_1,\ \ i=2,...,m.\tag
1.5.1$$ Thus, $$G_2=\frac{\sigma^y_m}{t^{m-1}}G_1,
G_3=\frac{\sigma^y_m}{t^{m-2}}G_2,...,
G_{i+1}=\frac{\sigma^y_m}{t^{m-i}}G_i, i=1,...,m-1.$$

 An elementary calculation shows that if we denote by
$V^m_i$ the 'mixed Van der Monde' matrix
$$V^m_i=\left[\matrix 1&\ldots&1\\x_1&\ldots&x_m\\\vdots&&\vdots\\
x_1^{m-i}&\ldots&x_m^{m-i}\\y_1&\ldots&y_m\\\vdots&&\vdots\\
y_1^{i-1}&\ldots&y_m^{i-1}\endmatrix\right]\tag 1.5.2$$then we
have
$$G_i=\pm\det (V^m_i).\tag 1.5.3$$
    In particular, $G_m$ as given in (1.5.1) coincides with
$v^m_y.$ I claim that $G_i$ generates $\O(-O\Gamma\spr m.)$ over
$\tilde{U}_i.$ This is clearly true where $t\neq 0$ and it remains
to check it along the special fibre $OH_{m,0}$ of $OH_m$ over $B$.
Note that $OH_{m,0}$ is a sum of components of the form
$$\Theta_I=\text{Zeros}(x_i,i\not\in I, y_i, i\in I),
 I\subseteq\{1,...,m\},$$ none of which is contained in the
 singular locus of $OH_m.$
  Set
 $$\Theta_i=\bigcup\limits_{|I|=i}\Theta_I.$$ Note that
 $$\Tilde{C}_i\times 0\subset \Theta_i, i=1,...,m-1$$ and therefore
 $$\tilde{U}_i\cap\Theta_j=\emptyset, j\neq i-1, i.$$
 Note that $y_i$ vanishes to order 1 (resp. 0) on $\Theta_I$
 whenever $i\in I$ (resp. $i\not\in I$). Similarly, $x_i-x_j$
 vanishes to order 1 (resp. 0) on $\Theta_I$ whenever both $i,j\in
 I^c$ (resp. not both $i,j\in I^c$). From this, an elementary
 calculation shows that the vanishing order of
 $G_j$ on every component $\Theta$ of
 $\Theta_k$ is $$\ord_{\Theta}(G_j)=(k-j)^2+(k-j).\tag 1.5.4$$
 We may unambiguously denote this number by $\ord_{\Theta_k}(G_j)$.
 Since this order is nonnegative for all $k,j,$
 it follows firstly that
 the rational function $G_j$ has no poles, hence is in fact regular
 on $X^m_B$ near
 $mp$ (recall that $X^m_B$ is normal); of course,
 regularity of $G_j$ is also immediate from (1.5.3).
  Secondly, since this order is zero for $k=j, j-1$, and
 $\Theta_j, \Theta_{j-1}$ contain all the components of $OH_{m,0}$
  meeting
 $\tilde{U}_j$, it follows that in $\tilde{U}_j,$ $G_j$ has no
 zeros besides $O\Gamma\spr m.\cap\tilde{U}_j,$ so $G_j$ is a generator
 of $\O(-O\Gamma\spr m.)$ over $\tilde{U}_j.$\par Now since $Z_j$ is a
 generator of $L$ over $\tilde{U}_j,$ we can define our
 isomorphism $\gamma$ over $\tilde{U}_j$ simply by specifying that
 $$\gamma (G_j)=Z_j\  {\text{on}}\ \tilde{U}_j.$$ Now to check that
 these maps are compatible, it suffices to check that
 $$G_j/G_k=Z_j/Z_k$$ as rational functions (in fact, units over
 $\tilde{U}_j\cap\tilde{U}_k$). But the ratios $Z_j/Z_k$ are
 determined by the relations (1.4.2-1.4.3), while $G_j/G_k$
 can be computed
 from (1.5.3), and it is trivial to check that these agree. This
 completes the proof of Theorem 1.\qed\par
 \proclaim{Corollary 4} The ideal of $OD^m$ is generated, locally
 near $p^m$, by $G_1,...,G_m.$\endproclaim\demo{proof} We have
 $$\I_{OD^m}=o\frak c_{m*}(\I_{O\Gamma\spr m.})=o\frak c_{m*}(L)$$
 is generated by
 the images of $Z_1,...,Z_m$, i.e.  by $G_1,...,G_m.$\enddemo
 As a further consequence, we can determine the ideal of the
 discriminant locus $D^m$ itself: let $\delta_m^x$
 denote the discriminant of $F_0$, which, as is
 well known \cite{L}, is a polynomial in the $\sigma_i^x$
 such that $$\delta_m^x=G_1^2.$$ Set
 $$\eta_{i,j}=\frac{(\sigma_m^y)^{i+j-2}}{t^{(i-1)(m-i)+(j-1)(m-j)}})
 \delta^m_x.\tag 1.5.5$$
 \proclaim{Corollary 5} The ideal of $D^m$ is generated, locally
 near $mp$, by $\eta_{i,j}, i,j =1,...,m.$\endproclaim\demo{proof} This
 follows from the fact that $\varpi_m$ is flat and that
 $$\varpi_m^*(\eta_{i,j})=G_iG_j, i,j=1,...,m$$ generate the ideal of
 $2OD^m=\varpi_m^*(D^m).$\enddemo
 Note that $\frak c_m^*(D^m)$ is a Cartier divisor
 on $X\sbr m._B$
 (that, of course, is just the universal property of blowing up) but
 its ideal, that is, $\O(-\frak c_m^*(D^m))$,
 is isomorphic in terms of our local model $\tilde{H}$ to
 $\O(2)$ (i.e. the pullback of $\O(2)$ from $\P^{m-1}$).
 This suggests that $\O(-\frak c_m^*(D^m))$
  is divisible by 2 as line bundle on $X\sbr m._B$,
 as the following result indeed shows. First some notation.
For a prime divisor $A$ on $X$, denote by $[m]_*(A)$ the prime
divisor on $X\sbr m._B$ consisting of schemes whose support meets
$A$. This operation is easily seen to be additive, hence can be
extended to arbitrary, not necessarily effective, divisors and
thence to line bundles.
 \proclaim{Corollary 6} Set $$\O_{X\sbr
 m._B}(1)=\omega_{X\sbr m._B}\otimes[m]_*(\omega_X\inv)\tag 1.5.6$$ Then

 $$\O_{X\sbr m._B}(-\frak c_m^*(D^m))\simeq\O_{X\sbr m._B}(2)\tag 1.5.7$$
 and
 $$\O_{X^{\lceil m\rceil}_B}(-o\frak c_m^*(OD^m))
 \simeq\varpi_m^*\O_{X\sbr m._B}(1).\tag 1.5.8$$

 \endproclaim

 \demo{proof} The Riemann-Hurwitz formula shows that the isomorphism
 (1.5.7) is valid on
 the open subset of $X\sbr m._B$ consisting of schemes disjoint from
 the locus of fibre nodes of $\pi$. Since this open is big (has
 complement of codimension $>1$), the iso holds on all of $X\sbr
 m._B$. A similar argument establishes (1.5.8)\qed\enddemo
 In practice, it is convenient to view (1.5.6)
 as a formula for $\omega_{X\sbr m._B},$
 with the understanding that $\O_{X\sbr m._B}(1)$ coincides in our
 local model with the $\O(1)$ from the $\P^{m-1}$ factor, and that
 it pulls back over $X^{\lceil m\rceil}_B=X\sbr m._B\times_{X\spr
 m._B}X^m_B$ to the $\O(1)$ associated to the blow up of the
 'half discriminant' $OD^m.$ We will also use the notation
 $$\O(\Gamma\spr m.)=\O_{X\sbr m._B}(-1),
 \Gamma\scl m.=\varpi_m^*(\Gamma\spr m.)$$
 with the understanding
 that $\Gamma\spr m.$ is Cartier, not necessarily effective, but $2\Gamma\spr
 m.$ and $\Gamma\scl m.$ are effective.
 \subheading{1.6 The small diagonal} The next Corollary
 will be crucial for the intersection calculus developed in the
 next section. It determines the restriction of the line
associated to $\Gamma\spr m.$, i.e. $\O_{X\sbr m._B}(1)$, on the
small diagonal.
 Thus let $\Gamma_{(m)}\subset
 X\sbr m._B$ be the small diagonal, which parametrizes schemes with
 1-point support, and which is the pullback of the small diagonal
  $$D_{(m)}\simeq X\subset X\spr m._B.$$ The
 restriction of the cycle map yields a birational morphism
 $$\frak c_m:\Gamma_{(m)}\to X$$ which is an isomorphism except
 over the set of fibre nodes sing$(\pi)$. Let $$J_m^\sigma\subset
 \O_X$$ be the ideal sheaf whose stalk at each fibre node is of
 type $J_m$ as in \S 0.

  \proclaim{Corollary 7} Via $\frak c_m$, $\Gamma_{(m)}$ is
  equivalent to the blow-up of $J_m^\sigma.$
  If $\O_{\Gamma_{(m)}}(1)$ denotes the canonical blowup
 polarization, we have
 $$\O_{\Gamma_{(m)}}(-\Gamma\spr m.)=\omega_{X/B}^{\otimes
 \binom{m}{2}}\otimes\O_{\Gamma_{(m)}}(1).\tag 1.6.1$$

 \endproclaim\demo{proof} We may work with the ordered versions
 of these objects, locally over a
 neighborhood of a point $p^m\in X^m_B$ where $p$ is a fibre node.
 There the ideal of $OD^m$ is generated by $G_1,...,G_m$ and $G_1$
 has the Van der Monde form $v^m_x$, while the other $G_i$ are given
 by (1.5.1). We try to restrict the ideal of $OD^m$ on the small
 diagonal $OD_{(m)}.$ To this end, note that
 $$(x_i-x_j)|_{OD_{(m)}}=dx=x\frac{dx}{x}$$
 and $\eta=\frac{dx}{x}$ is a local generator of $\omega_{X/B}.$
 Therefore
 $$G_1|_{OD_{(m)}}=x^{\binom{m}{2}}\eta^{\binom{m}{2}}.$$
 From (1.5.1) we then deduce
 $$G_i|_{\Gamma_{(m)}}=x^{\binom{m+i-1}{2}}y^{\binom{i}{2}}
 \eta^{\binom{m}{2}}, i=1,...,m.\tag 1.6.2$$
 Since $G_1,...,G_m$ generate the ideal  $I_{OD^m},$ it follows
 that
 $$I_{OD^m}\otimes\O_{OD_{(m)}}\simeq
 J_m^\sigma\otimes\omega^{\binom{m}{2}}.$$ Consequently, we also
 have $$I_{D^m}\otimes\O_{D_{(m)}}\simeq
 J_m^\sigma\otimes\omega^{\binom{m}{2}}.$$ Then pulling back to
 $X\sbr m._B$ we get (1.6.1).\qed\enddemo
 Now working locally at a point $p$ (which may be assumed a fibre
 node, though this is irrelevant for what follows), consider the
 blowup $c:\Gamma\to X$ of a punctual ideal of type $J_m$, and let $e_m$
 be the exceptional divisor, defined by
 $$\O_\Gamma(1):=\O_\Gamma(-e_m)=c^*J_m$$
 (pullback of ideals). Clearly the support of $e_m$ is
 $C^m=\bigcup\limits_{i=1}^{m-1}C^m_i,$ so we can write
 $$e_m=\sum\limits_{i-1}^{m-1}b_{m,i}C^m_i$$ and we have
 $$-e_m^2=\deg(\O(1).e_m)=\sum\limits_{i=1}^{m-1} b_{m,i}=:b_m.$$
 Now the general point on $C^m_i$ corresponds to an ideal
 $(x^{m-i+1}+ay^i), a\in\C^*$ and the rational function
 $x^{m-i+1}/y^i$ restricts to a coordinate on $C^m_i.$
 It follows that if
 $A_i\subset X$ is the curve with equation $f_i=x^{m-i+1}-ay^i, a\in\C^*,$
 then its
 proper transform $\tilde{A_i}$ meets $C^m$ transversely
 in the unique point
 $q\in C^m_i$ with coordinate $a$, so that
 $$\tilde{A_i}.e_m=b_{m,i}.$$ Thus, setting $J_{m,i}=J_m+(f_i)$
 we get following
 characterization of $b_{m,i}$:
 $$b_{m,i}=\l(\O_X/J_{m,i}).\tag 1.6.3$$
 To compute this, we start by noting that a cobasis $B_m$ for $J_m,$
 i.e. a basis for $\O_X/J_m$ is given by the monomials $x^ay^b$
 where $(a,b)$ is an interior point of the polygon $S_m$ as in \S 0;
 equivalently, the square with bottom left corner $(a,b)$ lies in
 $R_m.$ Then a cobasis $B_{m,i}$ for $J_{m,i}$ can be obtained by
 starting with $B_m$ and eliminating\par - all monomials $x^ay^b
 $ with $b\geq i$;\par - for any $j$ with $\binom{j}{2}\geq i,$
 all monomials that are multiples of
 $x^{\binom{m+1-j}{2}+m+1-i}y^{\binom{j}{2}-i}$;\nl the latter of
 course comes from the relation
 $$x^{\binom{m+1-j}{2}}y^{\binom{j}{2}}\equiv 0\mod J_m.$$
 Graphically, this cobasis corresponds exactly to the polygon
 $R_{m,i}$ in \S 0, hence $$b_{m,i}=\beta_{m,i},
 b_m=\beta_m;\tag 1.6.4$$
 in particular\proclaim{Corollary 8} With the above notations, we
 have globally $$e_m^2=-\sigma\beta_m,\tag 1.6.5$$
 $$\int\limits_{\Gamma_{(m)}}(\Gamma\spr
 m.)^2=-\sigma\beta_m+\binom{m}{2}^2\omega_{X/B}^2.\tag 1.6.6$$

\endproclaim
\remark{Remark 8.1} The components $C_{i}^m, i=1,...,m-1$ of $e_m$
 are special cases of the {\it{node scrolls}}, to be introduced in
 \S2.3 below; the general node scroll is a $\P^1$ bundle whose fibre
 is an $e_{m,i}.$ The coefficients $\beta_{m,i}$ play an essential role in the
 intersection calculus to be developed in \S2.\endremark
For the remainder of the paper, we set $$\omega=\omega_{X/B}\tag
1.6.7$$ (viewed mainly as divisor class).

 \heading {\bf{2. The tautological ring}}\endheading
 We continue with the notations and assumptions
of \S 1 and assume additionally that $X$ is a smooth surface and $B$
is a smooth curve. Our aim is to
 study the intersection theory associated to
 the tautological quotient bundle over the relative
 Hilbert scheme
 $X\sbr m._B$. Thus let $$\Lambda_m=\Spec
 (\O_{X\sbr m._B\times_BX}/\I_m)$$ be the universal length-$m$
 subscheme, and for any vector bundle $E$ on $X$, set
 $$\lambda_m(E)=p_{1*}(p_2^*(E)\otimes\O_{\Lambda_m}).$$
 By flatness of $\Lambda_m$ over $X\sbr m._B$,
 $\lambda_m(E)$ is clearly locally free of rank
 $m.\rk(E)$ on $X\sbr m._B$. Our plan is first to review a formula
 for (essentially)  the Chern
 classes of $\lambda_m(E)$, called {\it{tautological classes}}.
 More percisely,
 we will shift our situs operandi from the Hilbert scheme to its flag
 analogue. As a result, we are able to express the
 (pullback of the) tautological
 classes in terms of certain 'diagonal' {\it divisorial} classes (of
 Chow degree 1), essentially just the class $\Gamma\spr m.$
 defined above and its lower-degree analogues.
We then work out the products of tautological classes in the Chow
(or cohomology) ring of $X\sbr m._B,$ including  especially the
top-degree products, i.e.
 the Chern numbers of $\lambda_m(E)$, which might be
 called the {\it{tautological numbers}}. In the
 applications of the Hilbert scheme to classical enumerative geometry,
 it is these numbers that are required. We proceed, in fact,
by giving a set of {\it{additive generators}} for the ring generated
by the tautological classes $c_i(\lambda_m(E))$, and giving a
calculus for expressing the product of one of these generators by
$\Gamma\spr m.$ as a linear combination of other generators. This is
sufficient to compute all tautological numbers.
 \subheading{2.1 Divisorial multiplicative genrators}
 The total Chern
 class $c(\lambda_m(E))$ has been computed elsewhere in similar
 contexts: \cite{L} in the case of the (full) Hilbert scheme of a
 smooth surface, \cite{R} in the context of the relative
 flag-Hilbert scheme of a family of nodal curves over a base
 of any dimension.  %Here we will use the formula of \cite{R} to
 Our main goal is to
 compute the Chern numbers of $\lambda_m(E)$, and we note that
 Chern numbers, i.e. 'top' degree polynomials in the Chern classes,
 have a different meaning for the
 $(m+1)$-dimensional $X\sbr m._B$ than for the $2m-$dimensional
 Hilbert scheme of $X$. Accordingly Lehn's formula \cite{L}
 will be of no direct use to us. Rather, we will use
 the approach of \cite{R} which has the advantage of yielding degree-1
 (i.e. divisorial) multiplicative
 generators for the canonical ring, albeit at the cost of
 passing from the Hilbert scheme to its flag analogue.
  We now proceed to recall the required
 statement from \cite{R}.\par Let $$W^m=W^m(X/B)\overset\pi\spr m.\to\longrightarrow B$$ denote the relative
 flag-Hilbert scheme of $X/B$, parametrizing flags
 of subschemes $$z.=(z_1<...<z_m)$$
 where $z_i$ has length $i$ and $z_m$ is contained in some fibre
 of $X/B$. Let $$w^m:W^m\to X\sbr m._B, w\scl m.:W^m\to X\scl m._B$$ be the canonical
 (forgetful) maps. Let $$p_i:W^m\to X$$ be  the canonical map
 sending a flag $z.$ to the 1-point support of $z_i/z_{i-1}$
 and $$p^m=\prod p_i:W^m\to X^m_B$$ their (fibred) product,
 which might be called the 'ordered cycle map'.
   $W^m$ carries Cartier divisors
 $$\Delta\spr i.=\sum\limits_{j=1}^{i-1}\Delta^i_j$$ with each
 $\Delta^i_j$ a prime Weil divisor defined generically by
 $p_i(z.)=p_j(z.)$ (thus $\Delta\spr 1.=0$).
  We have $$w_m^*(\Gamma\spr m.)=(w\scl m.)^*(\Gamma\scl m.)=
  \sum\limits_{i=2}^m
 \Delta\spr i..\tag 2.1.1$$
  The formula of \cite{R}, Cor. 3.2 states that for any
 vector bundle $E$, we have $$c(w_m^*\lambda_m(E))=
 \prod\limits_{i=1}^m
 c(p_i^*(E)(-\Delta\spr i.))\tag 2.1.2$$ In particular, if $E=L$
 is a line
 bundle, we have
 $$c(w_m^*\lambda_m(L))=
 \prod\limits_{i=1}^m(1+[L\spr i.]-[\Delta\spr i.])\tag
 2.1.3$$ where $$L^{(i)}=p_i^*(L).$$ In \cite{R2}
 we showed that (2.1.3) can be
 used to derive a more 'explicit' sum-of-products
 formula for $c(\lambda_m(L))$
 on $X\sbr m._B$ which,
 when $X$ is a smooth surface,
  agrees with the
 restriction of a formula for the analogous bundles on
 $\Hilb_m(X)$ due to Lehn \cite{L}. For the purpose of computing
 Chern numbers, obviously either $W^m$ or $X\sbr m._B$ could be
 used since the set of numbers they yield differ by a factor of
 $m!$. We will work in the former context, where the simple product
 formula (2.1.2) holds. Note that this
formula has the added advantage of yielding directly the the
{\it{Chern roots}} of $w^*\lambda_m(L)$, which are useful in
computations.\par
  In view of (2.1.2), we call the subring $T^m=T^m(X/B)$
  of the $\Q$-Chow
 ring of $W^m$ generated by the $\Delta\spr i.$  and the
 $p_i^*(A^.(X)), i=2,...,m$ the {\it {tautological ring}} of $W^m$.
 In view of (2.1.1), we may replace the generators
 $\Delta\spr i., i=2,...,m$ by
 $\Gamma\spr i.$ or $\Gamma\scl i., i=2,...,m$
 which are more convenient (e.g.
 $\Gamma\spr i.$
 lives on $X\sbr i._B$). By their very definition, the various
 $T^m$'s form a chain $$T^2\to ...\to T^{m-1}\to T^m.$$
 Assuming
 $X$ is a surface, so that $\dim W^m=m+1,$ we will give a method,
 inductive in $m$, to express an arbitrary %degree-$(m+1)$
 nonzero monomial $M$ in $T^m$
 in terms of certain additive generators (to be specified below),
 assuming the analogous result in $T^{m-1}$ is known.
 %all similar integrals on $W^j, j<m$ are known,
 We may assume that $M$ is a monomial
 in $\Gamma\scl 2., ...,\Gamma\scl m.,$ hence expressible in
  the form $$M=M'(\Gamma\scl m.)^r$$
 with $M'\in T^{m-1}.$ By induction on $m$,  we may assume
 $M'$ is already expressed as a linear combination of
 the additive generators. Therefore, we may as well assume $M'$ is itself
 one of the additive generators in $T^{m-1}.$ Then, using induction on
 $r$,
  it will suffice to show
 how to express the product of an additive generator
 in $W^m$ by $\Gamma\scl m.$
 as a linear combination of additive generators.\par
 Now our additive generators for the tautological
 ring come in three flavors: the {\it{diagonals}}, analogous to
 Nakajima's creation operators; the {\it{node scrolls}},
 which are certain $\P^1-$bundles parametrizing schemes whose
 support contains some fibre
 nodes; and the {\it{node sections}}, which are certain cross-sections of
 node scrolls.
We first introduce the diagonal classes.
\subheading{2.2 Diagonal
classes}
 Note that for
 any pair of distinct pairs
 $(i<j)\neq (i'<j')$, the intersection $$\Delta_i^j\cap
 \Delta_{i'}^{j'}$$ is a well-defined codimension-2 cycle on $W^m$,
 because $\Delta_i^j$ and $
 \Delta_{i'}^{j'}$ are Cartier at the generic point of the
 intersection. Similarly, for any index-set $$I=(i_1<...<i_k)\subset [1,m]$$
 and any $c\in H^.(X),$ we have a well-defined cycle class that we call
 a {\it{connected diagonal monomial}} $$q_I[c]=c\spr
 i_1.\Delta_{i_1}^{i_2}\Delta_{i_2}^{i_3}...\Delta_{i_{k-1}}^{i_k}
 =c\spr i_1.\Delta_{I}.\tag
 2.2.1$$ When necessary to indicate the dependence on $m$ we'll sometimes
 write this as $q_I\spr m.[c]$. When $I$ is a singleton $\{i\}$, (2.2.1) reads
 $$q_i[c]=c\spr i..$$ $q_I[c]$ is an ordered analogue of Nakajima's
 creation operator $q_{|I|}[c]$ (cf. \cite{N, EG}).
 Likewise, for any partition $(I.)=(I_1,..., I_h)$, i.e.
 collection $I_1,...,I_h\subset [1,m]$ of
 pairwise disjoint subsets or 'blocks', with associated classes $c_1, ...,c_h$,
 we have a
 well-defined (disconnected, if $h>1$) {\it{diagonal monomial}}
 $$q_{(I.)}[c.]=q_{I_1}[c_1]\cdots q_{I_h}[c_h].$$ We view $(I.)$ as a sort of
 disconnected
 set with $I_1,...,I_h$ its connected components, and $(c.)$ as a locally constant
 $H^.(X)$-valued function on $(I.)$.. Note that
 $q_{(I.)}[c.]$ is supported on
 $$\Delta_{(I.)}=
 \Delta_{I_1}\cap\cdots\cap\Delta_{I_h}\sim q_{I_1}[1]\cdots q_{I_h}[1]$$
 which maps under the ordered cycle
 map to the appropriate diagonal locus $OD^m_{(I.)}$. It is
 obvious from (2.1.3) that the Chern classes of $w_m^*\lambda_m(L)$
are linear combinations of diagonal monomials. The coefficients are
worked out in \cite{R2}, and are consistent with Lehn's formula in
\cite{L}. We call the group generated by the diagonal monomials
$q_{(I.)}[(c.)]$ the group of {\it{diagonal classes.}}\par It is
worth noting that the diagonal classes $q_I[c]=q\spr m._I[c]$ behave
simply with respect to push-forward and pullback under the natural
map
$$\gamma^{m,m-1}:W^m\to W^{m-1}.$$ First,
$$I\subset [1,m-1]\Rightarrow
(\gamma^{m,m-1})^*q\spr m-1._I[c]=q\spr m._I[c]\tag 2.2.2$$
(consequently, it is safe to omit the superscript from $q\spr
m._I[c]$); next,
$$m\in I, |I|>1\Rightarrow
\gamma^{m.m-1}_*q_I[c]=q\spr m-1._{I\cap [1,m-1]}[c];\tag 2.2.3$$
$$I=(m)\Rightarrow\gamma^{m,m-1}_*(q_{(m)}[c])=\gamma^{m,m-1}_*(c\spr
m.)=(\pi\spr m-1.)^*\pi\spr m._*(c)\tag 2.2.4$$ (if $c$ is of Chow
degree 1 (cohomological degree 2), this is just $\deg_\pi(c)
1_{W^{m-1}}$ where $\deg_\pi(c)$ is the fibre degree).
$$\gamma^{m.m-1}_*q_I[c]=0, m\not\in I.\tag 2.2.5$$ Analogous
formulae hold also for diagonal monomials $q_{(I.)}[c.]$. By the
projection formula, it follows in particular that
$\gamma^{m.m-1}_*q_{(I.)}[c.]$ is a diagonal monomial in $T^{m-1}$.
Using this inductively, we see that that for any 0-dimensional
(degree-$(m+1)$) diagonal monomial $q_{(I.)}[c.]$, we can easily
compute the number $$\int\limits_{W^m}q_{(I.)}[c.].$$
 Unfortunately, the group generated by the diagonal classes is
not closed under multiplication by $\Delta^i$ or $\Gamma^i$ classes;
achieving closure requires introduction of node scroll and node
section classes.
 \subheading{2.3 Node scrolls} Consider a partition
 $$I_1\coprod I_2\coprod J_1...\coprod J_a\coprod
 K_1...\coprod K_b\subseteq [1,m]$$  such that
$$|I_1|,|I_2|>0$$ and that $I_1$ contains the smallest $I_1$
elements of $I_1\cup I_2$ in terms of the usual ordering on $[1,m];$
thus $I_1$ is an 'initial segment' of $I_1\cup I_2$.
  We will call
  $$\Phi=(I_1|I_2:J|K)$$
 a set of {\it{partition data}} with respect to $m$. More generally,
 if $I_1$ is not an initial segment of $I_1\cup I_2$, we identify
 $(I_1|I_2:J|K)$ with $(I_1'|I_2':J|K)$ where $I'_1$ is the initial
 segment of $I_1\cup I_2$ of cardinality $|I_1|$ and $I'_2=(I_1\cup
 I_2)\setminus I'_1.$
  The case where $J$ or $K$ is empty is
 included, and if both are empty we will write $\Phi=(I_1|I_2:)$.
 We think of $\Phi$ as indexing some of the variables $x_1, y_1,..., x_m,y_m$
 where $I_1, J$ (resp. $I_2, K$) refer to $y$ (resp. $x$) variables.
 \par $\Phi$ is said to be {\it full} if
 $$\bigcup\Phi :=I_1\cup...\cup K_b=[1,m].$$
 A {\it filling} of $\Phi$ is a full set of
 partition data $\Phi'=(I_1|I_2:J'|K')$ such that $J'$ (resp. $K'$)
 differs from $J$ (resp. $K$) only by 1-element blocks. We write
 this as $\Phi\prec\Phi'$
%Here $J=(J_1,...,J_a), K=(K_1,...,K_b)$ are pluri-multi-indices and
%we have used $\#I$ for $|I|$ for notational causes.
Let $X_1, ..., X_\sigma$ be the singular fibres of $\pi$. Assume
first that the singular fibre $X_s$ is a union of two smooth
components $X'_s, X^"_s$ meeting in a single point $n_s,$ and fix
$\Phi=(I_1|I_2:J|K)$. We also let $n_s', n_s"$ denote the preimages
of $n_s$ on $X'_s, X"_s,$ respectively. Set
$$X^\Phi_s=\prod\limits_{a}p_{J_a}\inv(\Delta_{X"_s})\times
\prod\limits_{a}p_{K_a}\inv(\Delta_{X'_s})\times\prod\limits_{i\in
I_1\cup I_2}p_i\inv(n_s)\tag 2.3.1$$ where $\Delta_{X'_s}\subset
(X'_s)^{K_a}$ etc. denotes the {\it small} diagonal. This depends on
$I_1, I_2$ only via $I_1\cup I_2$. Note that the irreducible
components of $X^\Phi$ are precisely the $X^{\Phi'}$ where $\Phi'$
is a filling of $\Phi$, and each of these is a smooth subvariety of
$X^m_B$, isomorphic to $(X_s")^{\l(J)}(X_s')^{\l(K)}$, where $\l(J)$
denotes the number of blocks (or 'length') of the partition $J$. Fix
such a filling $\Phi'$. We have
$$\sigma_i^y|_{X^{\Phi'}_s}=0, i>|J'|:=\sum |J'_j|,$$
$$\sigma_i^x|_{X^{\Phi'}_s}=0, i>|K'|.$$ In light of the relations
(1.4.1), which hold in a local model of the Hilbert scheme near the
'origin' $n_s^m$, this implies that the generic fibre of $p^m$ over
$X^{\Phi'}_s$ in a neighborhood of the 'origin', described in terms
of this model, has the form
$$\bigcup\limits_{i=|J'|+1}^{m-|I'|-1}C^m_i.\tag\dag$$
 Setting $r=m-|J'|-|K'|=|I_1|+|I_2|,$ this same
fibre can be identified, in terms of a local model of
 Hilb over a {\it generic} point of $X^{\Phi'}_s,$ with
 $$\bigcup\limits_{i=1}^{r-1}C^r_i\tag\ddag$$
 (cf. Remark 1.3.2).  We denote by
$$F_s^{\Phi'}\subset X\scl m._B\tag 2.3.2$$ the component of
$w\scl m.((p^m)\inv(X_s^{\Phi'}))$ with generic fibre
$C^m_{|J'|+|I_1|}$ in the first identification $(\dag)$ or
$C^r_{|I_1|}$ in the second $(\ddag)$. When there is no confusion,
we may use the same notation for the preimage of $F_s^\Phi$ in
$W^m,$ i.e. $(p^m)\inv(X_s^{\Phi'})$. This depends on $I_1, I_2$
only via $I_1\cup I_2, |I_1|$ (recall that $I_1$ is an 'initial
segment' of $I_1\cup I_2$). Informally, we think of $I_1$ and $J$
(resp. $I_2$ and $K$) as indexing $y$ (resp. $x$) variables, where
the $I$ variables are localized at the origin and the $J,K$
variables are free.\par The natural map
$$p^{\Phi'}:F^{\Phi'}_s\to X^{\Phi'}_s\tag 2.3.3$$ is a $\P^1$-bundle
(this is of course only true of the model of $F_s^\Phi$ in $X\scl
m._B$) .
 The locus $F^{\Phi'}_s$ is called the {\it{node
scroll}} corresponding to the node $s$ and the partition data
$\Phi'.$ We also set
$$F_s^\Phi=\bigcup\limits_{\Phi\prec\Phi'\ \text{full}}F_s^{\Phi'},
 \tag 2.3.4$$ \par We now
indicate the modifications needed to construct node scrolls for an
irreducible 1-nodal fibre $X_s.$ For a set of partition data
$\Phi=(I_1|I_2:J|K)$ we now insist that $K=\emptyset$ and that
$\Phi$ be full. Let $X_s'$ be the normalization of $X_s$, marked
with the 2 node preimage $n_s', n_s".$  Set
$$X_s^\Phi=\prod\limits_{a}p_{J_a}\inv(\Delta_{X'_s})
\times\prod\limits_{i\in I_1}p_i\inv(n"_s)\times \prod\limits_{i\in
I_2}p_i\inv(n'_s)\subset (X'_s)^m,
$$ which is the direct analogue of (2.3.1).
 Note that this locus has $2^{\l(J)}$ natural 'origins', viz. the
elements of
$$\prod\limits_{a}p_{J_a}\inv\{(n'_s)^{J_a},(n"_s)^{J_a}\}\times\prod\limits_{i\in
I_1}p_i\inv(n"_s)\times \prod\limits_{i\in I_2}p_i\inv(n'_s)
$$ where $(n'_s)^{J_a}$ is the diagonal point corresponding to $n'_s$ etc.
Let $$n:X_s^\Phi\to X_s^m\subset X^m_B$$ be the natural map induced
by normalization, and set $$F_s^\Phi=(p^m)\inv(n(X_s^\Phi)).$$ We
note that the restriction of $p^m$ lifts to a $\P^1$-bundle
projection
$$p^\Phi:F_s^\Phi\to X_s^\Phi.$$ Indeed, this may be checked locally
analytically on $X_s^\Phi$ and there is clear from our local
analytic model for the Hilbert scheme, in which the branches of
$X_s$ at the node already appear separated, and the target of the
cycle map appears as the product of the symmetric products of the
branches. Finally, set $$F^\Phi =\sum\limits_{s=1}^\sigma F_s^\Phi$$
(sum over all singular fibres, both reducible and irreducible).
\par

 A {\it{node
section}} class is by definition a class of the form $-\Gamma\scl
m..F^\Phi_s$. The group generated by the classes of node scrolls and
node sections is called the group of {\it{node classes}}. This group
and the operation of $\Gamma\scl m.$ on it will be sudied at length
in \S2.5.\par One obvious fact worth noting at the outset is that
for $\Phi=(I_1|I_2:J|K)$, if $i\in I_1\cup I_2$ and $c\in H^*(X)$ is
a class of positive degree (codimension), then
$$c\spr i..[F^\Phi_s]=0, \forall s$$ (e.g. because $c$ admits a
representative disjoint from sing$(\pi)$). It follows that
$$I\cap (I_1\cup I_2)\neq \emptyset, \deg(c)>0
\Rightarrow q_I[c].[F^\Phi_s]=0=q_I[c].(\Gamma\scl m..[F^\Phi_s]).
\tag 2.3.5$$

 \subheading{2.4 Cutting a diagonal class}
 Our aim in this subsection is to express the product of
  a diagonal class by $\Gamma\scl
m.$ as a linear combination of diagonal classes and node (scroll)
classes, generalizing the results of \S1.6. To this end, note that
for any multi-index (1-block partition) $I$ and any $i\in I$, the
projection
$$g=p_i:\Delta_I\to X$$ is independent of $i\in I$ and thus
$\Delta_I$ maps birationally, via the ordered cycle map, to
$$X\times_B(\prod\limits_{j\not\in I}\ _BX).$$ The generic fibre
of the induced map $\Delta_I\to\prod\limits_{j\not\in I}\ _BX$ is
isomorphic to the 'small diagonal' $\Gamma_{(|I|)}$ which
parametrized 1-point schemes. Recall that the intersection
$\Gamma\spr m..\Gamma_{(m)}=\Gamma\scl m..\Gamma_{(m)}$ was computed
in \S1.6. A similar reasoning shows that $\Delta_I.\Gamma\scl m.$
can be computed as the sum of the following terms
\par $\bullet\sum \Delta_{(I:(a,b))},\ \ \ $ the sum being over all $a<b$
with both $a,b\notin I$, where $(I:(a,b)) $ is the obvious 2-block
partition;\par $\bullet\sum\Delta_{I\cup\{a,b\}}, \ \ \ $ sum over
all $a<b$ with $|I\cap\{a,b\}|=1,$ where $I\cup\{a,b\}$ is the
obvious block;\par $\bullet q_I[\omega^\binom{|I|}{2}];\ \ \ $\par
$\bullet\sum\limits_{j=1}^{|I|-1}\beta_{|I|,j}
  F^{((i_1,...,i_j|i_{j+1},...,j_{|I|}:)}.$

In order to write this compactly, the following purely combinatorial
gadget will appear frequently below. Let $J=\{j_1, j_2\}, j_1\neq
j_2$ be an index pair, and $(I.)$ a partition. A new partition
$(I'.)=J\ltimes (I.)$ is obtained from $(I.)$ as follows.
\par - if $j_1\in I_a, j_2\in I_b$ for some $a,b$, remove $I_a,
I_b$ from $(I.)$ and inset $I_a\cup I_b$ (in other words, 'connect
up' $I_a$ and $I_b$, reducing the number of blocks (or 'connected
components') by 1;\par - if $j_1\in I_a, j_2\notin I_b, \forall b$,
or vice versa, replace $I_a$ by $I_a\cup J$;
\par -
 if $j_1, j_2\notin I_a, \forall a,$ insert $J$ to
$(I.)$ as a block (thus increasing by 1 the number of connected
components).
\par - if $J\subset I_a$ for some $a$, $(I'.)=(I.).$\par
 With this notation, we can rewrite our formula for $\Gamma\scl
 m..\Delta_I$ as follows.
 $$\Gamma\scl m..\Delta_I=\sum\limits_{ i<j}
 \Delta_{\{i,j\}\ltimes I}+\sum\limits_{j=1}^{|I|-1}\beta_{|I|,j}
  F^{((i_1,...,i_j|i_{j+1},...,j_{|I|}:)}
  +\binom{|I|)}{2}q_I[\omega].\tag 2.4.1$$ The extension of (2.4.1)
  to the case of (disconnected) diagonal monomials is straightforward.
  For notational economy it is convenient to denote the middle term
  in (2.4.1)
  by $F^{(I:)}$; we similarly have $F^{(I:J|K)}$
  for partitions $(I:J|K).$
 From this it is easy to see that
more generally, we have
 $$\Gamma\scl m..\Delta_{(I.)}=\sum\limits_{i<j}
 \Delta_{\{i,j\}\ltimes(I.)}
 +\sum\limits_k\sum\limits_{J\cup K=I.\setminus I_k}F^{(I_k:J|K)}
 \tag 2.4.2$$
 $$+\sum\limits_kq_{I_1}[1]\cdots
 q_{I_k}[\binom{|I_k|}{2}\omega]\cdots q_{I_h}[1]$$
 where the last two sums may be
 restricted to those $I_k$ such that $|I_k|\geq 2,$ as the
 others yield 0 and, as always, for an irreducible singular fibre
 $X_s$ the condition that $K=\emptyset$ in $F_s^\Phi$ remains in force.
 For instance, in terms of the generator $G_1,$
the first term in (2.4.2) corresponds to the factors $x_i-x_j$ of
$G_1$ such that $\{i,j\}$ are not in the same block of $(I.)$; the
second and 3rd terms come from the various $k$ so that
$\{i,j\}\subset I_k.$
 \par It is a routine matter, albeit necessary,
 to extend (2.4.2) to a formula for
 diagonal monomials $\Gamma\scl m..q_{(I.)}[c.]$. To state this, we need yet some more
 notation. For any pluri-multi-index $(I.)=(I_1,...,I_h)$
 and classes $c_1,...,c_h\in H^*(X)$ (or $A^*(X)$), let us denote the diagonal
 monomial $q_{I_1}[c_1]
 \cdots q_{I_h}[c_h]$ by $q_{(I.)}[(c.)]$. Here the pluri-class
 $(c.)$
 should be viewed as a function from $(I.)$ to $H^*(X).$
  Then for a distinct pair $J=\{j_1, j_2\}$,
 there is a natural way to modify $(c.)$ to
  define a pluri-class $J\ltimes (c.)$ on
 $J\ltimes (I.)$:\par - in case $I_a$ and $I_b$ get connected up
 to form $I_a\cup I_b$, i.e. $j_1\in I_a, j_2\in I_b$
 or vice versa,  the value of $J\ltimes (c.)$ on $I_a\cup
 I_b$ is $c_a._Xc_b$;\par - in case $I_a$ gets
 replaced by $I_a\cup J$, the value of $J\ltimes(c.)$ on
$I_a\cup J$ is $(c_a)$;\par -
   in case $J$ is inserted to
 $(I.)$, define the value of $J\ltimes(c.)$ on $J$ to equal
 $1\in H^*(X)$;\par
 all other values are carried over from $(c.)$ to $J\ltimes(c.)$ in
 the obvious way.\par Also, if $(c.)$ is a pluri-class on $(I:J|K)$,
 define $F_s^{(I:J|K)}[(c.)]$ as follows.
  $$F_s^{(I:J|K)}[(c.)]=0\ \ \text{if}\ \ \deg c(I)>0;$$
$$F_s^{(I:J|K)}[(c.)]=F^{(I:J|K)}\prod c(J_a)\spr \min(J_a).
\prod c(K_a)\spr\min(K_a).\ \  \text{if}\ \  c(I)=1.\tag 2.4.3$$
 These are called generalized node scroll
classes, and we similarly have generalized node section classes.
Note that (2.4.3) clearly vanishes if $c(J_a)$ or $c(K_a)$ is of
degree $>1.$
 Also set $$X_s^{(I:J|K)}[(c.)]=X_s^{(I:J|K)}\prod
c(J_a)\spr \min(J_a). \prod c(K_a)\spr\min(K_a)..$$ Note that this
is a 0-cycle precisely when
$$\l(J)+\l(K)=\sum\limits_a\deg(c(J_a))+\sum\limits_a\deg(c(K_a))$$
where $\l(J), \l(K)$ denote the number of blocks in the partition
(which coincides with the dimension of $X_s^{(I:J|K)}$); in other
words, $X_s^{(I:J|K)}[(c.)]$ is a 0-cycle precisely when each
$c(J_a), c(K_a)$ is of degree 1. In this case we have
$$\int_{W^m}X_s^{(I:J|K)}[(c.)]=\int_{X_s^{(I:J|K)}}\prod
c(J_a)\spr \min(J_a). \prod
c(K_a)\spr\min(K_a).$$$$
=\prod\limits_a\deg_\pi(c(J_a))\prod\limits_a\deg_\pi(c(K_a))
 $$\nl With this notation, the extension of (2.4.2) reads

 $$\Gamma\scl m..q_{(I.)}[c.]=\sum\limits_{1\leq i<j\leq
 r}q_{\{i,j\}\ltimes(I.)}[\{i,j\}\ltimes(c.)]\tag 2.4.4$$
 $$+\sum\limits_k\sum\limits_{J\cup K=I.\setminus I_k}F^{(I:J|K)}
 [(c.)]+\sum\limits_kq_{I_1}[c_1]\cdots
 q_{I_k}[\binom{|I_k|}{2}\omega.c_k]\cdots q_{I_h}[c_h]
 $$ We have thus shown that the product of any diagonal class with
 $\Gamma\spr m.$ can be expressed in terms of diagonal classes and
 (generalized) node classes.
 \comment
 Using the relations ?, one can check that ? holds over the other
 open sets $U_j$ where $G_j$ is a generator. Consequently, we have
 $$\Gamma\spr r..q_{(I.)}[(c.)]=\sum\limits_{1\leq i<j\leq r}$$
\endcomment
We can now state the main result of this section:\proclaim{Theorem
4} Any element of the tautological ring $T^m$ can be (computably)
expressed as a linear combination of diagonals and  generalized node
classes.\endproclaim The plan is to prove by induction on $m$, so we
may assume it holds for all $m'<m.$ Note that the minimum dimension
for a generalized node scroll class $F_s^\Phi[(c.)]$ (resp.
generalized node section $-\Gamma\scl m..F_s^\Phi[(c.)]$) is 1
(resp. 0), both achieved when $X^\Phi_s.[(c.)]$ is 0-dimensional, so
in view of the obvious fact, when $X^\Phi_s([(c.)]$ is a 0-cycle,
that
$$\int_{W^m}-\Gamma\scl
m..F_s^\Phi[(c.)]=\int\limits_{X^m_B}X^\Phi_s([(c.)])\tag 2.4.5$$
(the latter being the degree of a 0-cycle)  Theorem 4 allows us to
compute $\int_{W^m}M$ for any top-degree element $M\in T^m$, as was
our main goal. \subheading{2.5 Cutting a node class} It remains to
analyze the product of a generalized node class with $\Gamma\spr m.$
(i.e. with $\Gamma\scl m.$). We will do this for  ungeneralized node
classes, as the extension to the case of generalized node classes is
straightforward.
 To this end, we wish first to
analyze the structure of a node scroll $F_s^\Phi$ with
$\Phi=(I.:J|K)$ a full set of partition data. To be able to state
formulae uniformly the reducible and irreducible singular fibres, it
is convenient to set
$$\matrix K'&=K,\  \text {reducible case}\\&=J,\  \text{irreducible
case}\endmatrix$$  $$\matrix K"&=K,\  \text {reducible
case}\\&=\emptyset,\ \text{irreducible case}\endmatrix$$As noted
earlier, the natural map
$$p^\Phi:F_s^\Phi\to X_s^\Phi$$ exhibits $F_s^\Phi$ as a $\P^1$-bundle,
and we wish to identify the corresponding vector bundle. Assume to
simplify notation that $I_1=[1,i], I_2=[i+1, r].$ Recall that
homogeneous coordinates on $C^r_i$ are given by $Z_i, Z_{i+1}$ which
correspond to the mixed Van der Monde generators $G_i, G_{i+1};$
ditto for $C^m_{|J|+i}.$ Consider the mixed Van der Monde matrix
$V^m_{|J|+i}$ whose determinant yields $G_{|J|+i}.$ It has an
$r\times r$ block submatrix based on the $I$-indexed rows and the
columns, corresponding to $1, x, ..., x^{r-i}, y,..., y^{i-1}$,
whose determinant is equal to $G_{i,I},$ that is, the $G_i$
expression in the variables $x_1, y_1,..., x_r, y_r.$ Note that this
is globally defined along $X_s^\Phi$. The determinant of the
complementary submatrix, considered as function on $X_s^\Phi$, is a
'shift' of another Van der Monde, equal to
$$(x^{K'})^{r-i}(y^J)^{i}\prod\limits_{a<b\in
\bigcup K'}(x_a-x_b)\prod\limits_{a<b\in \bigcup J}(y_a-y_b),\tag
2.5.1$$ where $x^{K'}=\prod\limits_{k\in \bigcup K'}x_k$ etc and,
for $X_s$ irreducible, $x, y$ are local coordinates at the node
preimages $n'_s, n"_s$, respectively. Note that in the irreducible
nodal case,  the last 2 factors in (2.5.1) define the same diagonal
locus, the one near $(n'_s)^{K'}$, the other near $(n"_s)^J$. Now
(2.5.1) is a generator of the invertible ideal
$$'E_s^\Phi=\O(-(r-i)\sum\limits_{a\in \bigcup K'}
p_a^*n'_s-i\sum\limits_{a\in \bigcup J} p_a^* n"_s-\sum\limits_{a,b
\in \bigcup K"} p_{a,b}^*(\Delta)-\sum\limits_{a,b \in \bigcup J}
p_{a,b}^*(\Delta)).\tag 2.5.2$$ Other terms in the Laplace expansion
of $G_{|J|+i}$ along the $I$ columns have order
$>\binom{r}{2}=\ord(G_{i,I})$ in the $I$ variables. Analogous
considerations for the second Van der Monde generator $G_{|J|+i+1}$
lead to the invertible ideal
$$"E_s^\Phi=\O(-(r-i-1)\sum\limits_{a\in \bigcup K'}
p_a^*n'_s-(i+1)\sum\limits_{a\in \bigcup J} p_a^*
n"_s-\sum\limits_{a,b \in \bigcup K"}
p_{a,b}^*(\Delta)-\sum\limits_{a,b \in \bigcup J}
p_{a,b}^*(\Delta)).\tag 2.5.3$$ Setting $$E^\Phi_s= 'E^\Phi_s\oplus
"E^\Phi_s,\tag 2.5.4$$ we conclude that, at least in a neighborhood
of the 'origin' $n_s^m$, we have
$$F^\Phi_s\simeq\P(E^\Phi_s)\tag 2.5.5$$ so that
$$\O(-\Gamma^m)|_{F^\Phi_s}=\O_{\P(E^\Phi_s)}(1).\tag 2.5.5'$$ A similar
argument shows that this isomorphism persists near 'less special'
points on $X^\Phi_s$, namely, expanding $G_{|J|+i}$ we again get,
modulo higher-order terms,  the same $G_{i,I}$ factor times another
local generator of $'E^\Phi_s$ and likewise for $G_{|J|+i+1}$; so
the isomorphism (2.5.5)-(2.5.5') holds globally. Note that
$$\P(E^\Phi_s)=\P(\O(-\sum\limits_{a\in J}
p_a^*n"_s)\oplus\O(-\sum\limits_{a\in K'} p_a^* n'_s))\tag 2.5.6$$
but the latter bundle gives the 'wrong' $\O(1)$.\par Next it is
important to compare node classes on $W^{m-1}$ and $W^m$. Let
$\Phi=(I_1|I_2:J|K)$ be full partition data with respect to $[1,
m-1].$ In the reducible case, there are precisely two completions of
$\Phi$ with respect to $[1,m]$, namely
$$\Phi'=(I_1|I_2:J^+=J\cup\{m\}|K), \Phi"=(I_1|I_2:J|K^+=K\cup\{m\}).$$
In the irreducible case, there is just $\Phi'$. There is a natural
sheaf inclusion
$$E_s^{\Phi'}\to p_{[1,m-1]}^*E_s^\Phi(-ip_m^*(n_s)-\sum\limits_{a\in\bigcup
J}p^*_{a,m}(\Delta))\tag 2.5.7$$ which drops rank by 1 with
multiplicity 1 along $p_m\inv(n_s)$, identifying $F_s^{\Phi'}$ as an
elementary modification of $F_s^\Phi$, albeit with polarization
$$-\Gamma\spr m..F_s^{\Phi'}=
(-\Gamma\spr m-1.-(i+1)p_m^*(n_s)-\sum\limits_{a\in\bigcup
J}p^*_{a,m}(\Delta)).F_s^\Phi\tag 2.5.8$$ (see Remark 2.5.1 below).
In the reducible case, we have additionally $$-\Gamma\spr
m..F_s^{\Phi"}= (-\Gamma\spr
m-1.-(i+1)p_m^*(n_s)-\sum\limits_{a\in\bigcup
K}p^*_{a,m}(\Delta)).F_s^\Phi.\tag 2.5.8"$$ In fact the model of
$F_s^{\Phi'}$ on $W^m$ is a blown-up $\P^1$-bundle which contracts
on the one hand to $F_s^\Phi\subset X\scl m-1._B$ and on the other
hand to $F_s^{\Phi'}\subset X\scl m._B$. Together with (2.5.8) and
(2.5.8"), this implies that $(\gamma^{m,m-1})^*$ takes node classes
on $W^{m-1}$ to node classes on $W^m.$ From this, it is obvious that
the same is true for generalized node classes. Now to compute the
Chern classes of $E^\Phi_s$, note that
$$p_a\inv(n"_s)=[X_s^{(I_1\cup
\{a\}|I_2:J\setminus\{a\}|K)}], a\in\bigcup J \tag 2.5.9$$
$$p_a\inv(n'_s)=[X_s^{(I_1|I_2\cup
\{a\}:J|K\setminus\{a\})}], a\in\bigcup K'\tag 2.5.10$$ in the
reducible case, and $$p_a\inv(n'_s)=[X_s^{(I_1|I_2\cup
\{a\}:J\setminus\{a\}|K)}], a\in\bigcup K'\tag 2.5.10'$$ in the
irreducible case;

$$p_{a,b}^*(\Delta)=\omega\spr
a.=\omega\spr\min(J_r).=(2g(X_s")-2)p_a^*(pt) \ \text{if}\
\{a,b\}\subset J_r\tag 2.5.11$$
$$p_{a,b}^*(\Delta)=\omega\spr a.=\omega\spr\min(K'_r).=(2g(X_s')-2)p_a^*(pt)
\ \text{if}\ \{a,b\}\subset K'_r.\tag 2.5.12$$
$$p_{a,b}^*(\Delta)=
[X_s^{(I.:(a,b)\ltimes J|K)}]\  \text{if}\ \{a,b\}\not\subset
J_r,\forall r, \{a,b\}\subset\bigcup J\tag 2.5.13$$
$$p_{a,b}^*(\Delta)=
[X_s^{(I.:J|(a,b)\ltimes K")}] \ \text{if}\ \{a,b\}\not\subset K_r,
\forall r, \{a,b\}\subset \bigcup K"\tag 2.5.14$$(here just (2.5.12)
and (2.5.13) are operative in the irreducible case).\nl All these
are codimension-1 classes on $X^\Phi_s$, whose pullback via $p^\Phi$
are clearly generalized node classes. It follows that, in the
irreducible case, \nl$c_1(E^\Phi_s)=$$$
-(2i+1)(p^\Phi)^*\sum\limits_{a\in J}[X_s^{(I_1\cup
\{a\}|I_2:J\setminus\{a\}|K)}]-(2r-2i-1)(p^\Phi)^*\sum\limits_{a\in
K}[X_s^{(I_1|I_2\cup\{a\}:J|K\setminus\{a\})}] $$
$$-2\sum\limits_{a<b\in\bigcup K}(p^\Phi)^*[X_s^{(I.:J|(a,b)\ltimes K)}]
-2(p^\Phi)^*[X_s^\Phi]\sum\limits_r\binom{|K_r|}{2}\omega\spr
\min(K_r).
$$$$-2(p^\Phi)^*\sum\limits_{a<b\in\bigcup J}[X_s^{(I.:(a,b)\ltimes J|K)}]
-2(p^\Phi)^*[X_s^\Phi]\sum\limits_r\binom{|J_r|}{2}\omega\spr
\min(J_r).
$$

$$= -(2i+1)\sum\limits_{a\in J}[F_s^{(I_1\cup
\{a\}|I_2:J\setminus\{a\}|K)}]-(2r-2i-1)\sum\limits_{a\in
K}[F_s^{(I_1|I_2\cup\{a\}:J|K\setminus\{a\})}] $$
$$-2\sum\limits_{a<b\in\bigcup K"}[F_s^{(I.:J|(a,b)\ltimes K")}]
-2[F_s^\Phi]\sum\limits_r\binom{|K_r|}{2}\omega\spr \min(K"_r).
$$$$-2\sum\limits_{a<b\in\bigcup J}[F_s^{(I.:(a,b)\ltimes J|K)}]
-2[F_s^\Phi]\sum\limits_r\binom{|J_r|}{2}\omega\spr \min(J_r).; \tag
2.5.15 $$ in the irreducible case, the second summation is replaced
by
$$\sum\limits_{a\in
K'}[F_s^{(I_1|I_2\cup\{a\}:K'\setminus\{a\}|\emptyset)}] $$ In the
expression (2.5.15) the 1st 2 terms come from the 1st 2 terms in
$'E, "E$; the 3rd and 4th terms come from the 3rd term in $'E, "E$
and correspond to the case where $a,b$ are in different blocks
(resp. the same block) of $K$; similarly for the 5th and 6th terms.
In particular, $c_1(E^\Phi_s)$ is clearly a generalized node class.
The computation of
$$c_2(E^\Phi_s)=c_1('E^\Phi_s)c_1("E^\Phi_s)\tag 2.5.16$$ is straightforward:
note that $X_s^\Phi$ is just a product of smooth curves and the
classes being multiplied are standard ones. The following elementary
facts may be used:
$$p_a^*(n_s)^2=0;\tag 2.5.17$$
$$p_a^*(n_s)p_b^*(n_s)=p_{a,b}^*(pt)=p_{a,b}^*(\Delta)p_a^*(n_s),
a\neq b$$$$p_a^*(n"_s)p_b^*(n"_s)=X^{(I_1\cup\{a,b\}|I_2:J|K)},
\{a,b\}\subset\bigcup J
$$
$$p_a^*(n'_s)p_b^*(n'_s)=X^{(I_1|I_2\cup\{a,b\}:J|K)}, \{a,b\}\subset\bigcup K'\tag 2.5.18
$$$$p_a^*(n"_s)p_b^*(n'_s)=X^{(I_1\cup\{a\}|I_2\cup\{b\}:J|K)}, a\in\bigcup J,b\in\bigcup
K';
$$
$$p_{a,b}^*(\Delta)p_{c,d}^*(\Delta)=0;\tag2.5.19
$$ if $a,b,c,d$ are in the same block of $J$ or $K$;\nl

if $a,b$ are in different blocks, then
$$
p_{a,b}^*(\Delta)^2=(2-2g)p_{a,b}^*(pt);\tag 2.5.20$$ where
$g=g(X")$ if $a,b\in\bigcup J$ or $g(X')$ if $a,b\in\bigcup K';$\nl
more generally, if $a,b$ are in different blocks of $J$, then for
all $c,d$,
$$p_{a,b}^*(\Delta)p_{c,d}^*(\Delta)=p_{c,d}^*(\Delta)|_{X^{(I_1|I_2:(a,b)\ltimes
J|K)}};\tag 2.5.21$$ ditto if $a,b$ are in different blocks of $K$.

 Clearly $c_2(E^\Phi_s)$ is in the group of generalized
node classes. Now  Grothendieck's standard relation

$$c_2(E_s^\Phi(-1))=0$$ yields
$$(\Gamma\spr m.)^2.F_s^\Phi=-\Gamma\spr
m..(p^\Phi)^*c_1(E^\Phi_s)-(p^\Phi)^*c_2(E^\Phi_s).\tag 2.5.22$$
Therefore also
$$(\Gamma\spr m.)^2.F_s^\Phi[(c.)]=\tag 2.5.23$$$$-\Gamma\spr
m..(p^\Phi)^*(c_1(E^\Phi_s).[X^\Phi_s[(c.)]-(p^\Phi)^*(c_2(E^\Phi_s).[X^\Phi_s[(c.)]
.$$ Applying this recursively, we see that the group of generalized
node classes is closed under multiplication by $\Gamma\spr m.$,
which completes the proof of Theorem 4.\qed

\comment any class $\delta=(\Gamma\spr m.)^r.F_s^\Phi$ for any $r$
can be computed as $\Gamma\spr m..(p^\Phi)^*\alpha+(p^\Phi)^*\beta$
where $\alpha, \beta$ are node classes. Continuing in this process
until $\delta$ is 0-dimensional, we must have in that case that
$\alpha$ is zero-dimensional and $\beta=0$ so the evident fact
$$-\Gamma\spr m..(p^\Phi)^*(pt.)=1$$ concludes the computation.

\par

ZZZZZZZZZZZZZZZZZZZZZZZZZZZZZZZZZZZZZZZZZZ

 Write $M$ in the form
 $$M=A_1...A_m,$$ $$A_i=(c_i)^{(i)}(\Delta^{(i)})^{r_i}, c_i\in
 A^.(X), r_i\geq 0.$$ \proclaim{Rule 1} For each
 $i=2,...,m,\ $ $A_i$ has positive degree.\endproclaim
 Indeed
$$\int\limits_{W^m}M=\int\limits_{W^{i-1}}\gamma^{i,i-1}_*\gamma^{m,i}_*(M)$$
and if $A_i=1$, then the
 projection of $\gamma^{m,i}_*(M)$ to $W^{i-1}$, as cycle, has 1-dimensional
 fibres, hence the integral vanishes.\proclaim{Rule 2} With at
 most 2 exceptions, each $A_i, i=1,...,m$ has degree 1.\endproclaim
  This
 is immediate from Rule 1.\proclaim{Pre-rule 3}
 If $r_m=0,$ then $\int\limits_{W^m} M$ is computable from
 integrals of tautological ring elements on
 $W^{m-1}.$\endproclaim Indeed  in this case $\deg(c_m)=1$ or $2$ .
 If $\deg(c_m)=1$, then
 $$\int\limits_{W^m}M=
 \deg_\pi(c_m)\int\limits_{W^{m-1}}A_1\cdots A_{m-1}$$ where
 $\deg_\pi(c):=c.f, f:=\pi\inv(pt.).$ If $\deg(c_m)=2, $ so $c_m=k[pt]$,
 then
 $$M=k\int\limits_{W^{m-1}}A_1\cdots A_{m-1}f\spr m-1.
 .$$
 \proclaim{Rule 3} If $r_i=0$ for any $i\in[2,m]$,
 then $\int\limits_{W^m} M$ is computable from
 integrals of tautological ring elements on
 $W^{m-1}.$\endproclaim
  Let $i$
 be largest with $r_i=0$. We may assume $i<m.$ Since $c_i$ may be
 represented by a cycle disjoint from sing$(\pi)$, the class
 $c_i\spr i.\Delta\spr j._i, $ for any $j>i$,  is well defined.
 Let us write, in the
 obvious sense, $$\Delta\spr j.=\Delta\spr j._i+\Delta\spr
 j._{\widehat{i}}.$$
 %Assume first that $r_{i+1}=1.$
 Then we can
 rewrite $M$ as a sum of monomials
 $$M=M_0+M_1+...$$
 where $M_0$ involves only $\Delta\spr j._{\widehat{i}}$ for various $j$
  and no $\Delta\spr
 j._{{i}}$ for any $j$,  and the other monomials involve some $\Delta\spr
 j._{{i}}.$ Now $M_0$ may be integrated by a simple shuffle:
 replace each $j=i+1,...,m$ by $j-1$ and $i$ by $m$. This
 transforms $\Delta\spr j._{\widehat{i}}$ into $\Delta\spr j-1.$,
 thus yielding a 'standard'
 monomial $M'_0$ which may be integrated by Pre-rule 3
 above.\par To integrate the other monomials, pick one, say $M_1,$
 and let $j$ be smallest so that $\Delta\spr j._i$ occurs in
 $M_1.$
 Suppose first that
  $\Delta\spr j._i$ appears in $M_1$ precisely to 1st
 power. Then we can transform $M_1$ as follows.
 Write $$c_i\spr i.\Delta\spr
 j._i=c_i\spr j.,$$ then replace all occurrences of $i$ in $M_1$
 (e.g. in a factor $\Delta\spr k._i$)  by
 $j$, then concatenate by replacing each $k>i$ by $k-1$. This
 yields a monomial $M'_1$ of degree $m$ in $W^{m-1}$ such that
 $$\int\limits_{W^m}M_1=\int\limits_{W^{m-1}}M'_1\tag 2.3$$
and the latter is computable by
 induction.\par Next, $\Delta\spr j._i$ appears in $M_1$ precisely to 2nd
 power, for some $j$, we can write $$c_i\spr i.(\Delta\spr
 j._i)^2=(c_i.\omega)\spr j., \omega=\omega_{X/B}$$ (again
 because $c_i$ may be assumed supported away from sing$(\pi)$),
  and again
 replace all occurrences of $i$ in $M_1$ by
 $j$, then concatenate by replacing each $k>i$ by $k-1$, again
 obtaining a monomial $M'_1$ on $W^{m-1}$ for which the analogous
 exponent $r'_m=0$ holds and which may be integrated via pre-rule 3.
 \par Finally if
$\Delta\spr j._i$ appears in $M_1$ to $>2$nd
 power, for some $j$, we can write
$$c_i\spr i.(\Delta\spr
 j._i)^3=(c_i.\omega^2)\spr j.=0$$
 for dimension reasons, since $\deg (c_i)>0$ by assumption; therefore
 $\int\limits_{W^m}M_1=0$ in this case.
 We have verified Rule 3.\qed\par In light of Rules 1-3,  we
 may assume $M$ is divisible by every $\Delta\spr i., i=2,...,m$,
 hence $M$ is (rationally) divisible by the small diagonal
 $$\Delta_{(m)}=\frac{1}{(m-1)!}
 \Delta\spr 2.\cdots\Delta\spr m.$$of degree $m-1$.
  Write
 $$M=\Delta_{(m)}.B_1B_2$$ where each of $B_1, B_2$ is either a
 $\Delta\spr i.$ or a $c\spr i.$. There are thus 3 cases, and we
 begin with the more difficult one.
%\comment
 Note that

 $$\Delta_{(m)}=\frac{1}{m!}w^*(\Gamma_{(m)}),$$ hence
 $$\int\limits_{W^m}M=\frac{1}{m!}\int\limits_{X\sbr m._B}
 \Gamma_{(m)}.w_*(B_1)w_*(B_2).$$

 \proclaim{Case 1}
  $$B_1=\Delta\spr i.,
 B_2=\Delta\spr j., i\leq j$$\endproclaim
\noindent{\bf{Case 1a:}} $j<m.$\nl Then we have
$$\int\limits_{W^m}M=\int\limits_{W^{m-1}}\Delta_{(m-1)}.
\Delta\spr i..\Delta\spr j.\tag 2.4$$ reducing the integral to one
for lower $m$.\nl {\bf{Case 1b:}}  $i<j=m.$\par Then note the
'telescope identity'
$$\gamma^{m,m-1}_*(\Delta\spr m.)^2=-\sum\limits_{k=1}^{m-1}\omega\spr k.
+2\sum\limits_{k=2}^{m-1}\Delta\spr k.\tag 2.5$$ (compare
\cite{R}, Remark 2.7.1). Indeed since (2.5) is an equality of
divisors, it suffices to check it off the codimension-2
exceptional locus of the ordered cycle map $$p^{m-1}:W^{m-1}\to
X^{m-1}_B,$$ where it is obvious. From (2.5) we get
$$\int\limits_{W^m}M=\int\limits_{W^m}\Delta_{(m-1)}\Delta\spr i.
(\Delta\spr m.)^2=$$$$=
\int\limits_{W^{m-1}}\Delta_{(m-1)}.\Delta\spr i.
(-\sum\limits_{k=1}^{m-1}\omega\spr k.
+2\sum\limits_{k=2}^{m-1}\Delta\spr k.) \tag 2.6$$ again reducing
the integral to a lower-dimensional one.\nl {\bf{ Case 1c:}}
i=j=m.\par Note that $$w^*\Gamma\spr
m.=2\sum\limits_{i=2}^m\Delta\spr i.,$$
$$(w^{m-1})^*(\gamma^{m,m-1})^*\Gamma\spr
m-1.=2\sum\limits_{i=2}^{m-1}\Delta\spr i.,$$ where
$w^{m-1}:W^{m-1}\to X\sbr m-1._B$ is the natural map. Abusing
notation, we will write these divisors respectively as $\Gamma\spr
m., \Gamma\spr m-1.$.

 Then we can write
$$4\int\limits_{W^m}M=4\int\limits_{W^m}\Delta_{(m)}(\Delta\spr
m.)^2= \int\limits_{W^m}\Delta_{(m)}(\Gamma\spr m.-\Gamma\spr
m-1.)^2$$$$=\int\limits_{W^m}\Delta_{(m)}(\Gamma\spr
m.)^2-2\int\limits_{W^m}\Delta_{(m)}\Gamma\spr m.\Gamma\spr m-1.
+\int\limits_{W^m}\Delta_{(m)}(\Gamma\spr m-1.)^2.$$$$=I+II+III$$
The second term II, writing $\Gamma\spr m.=\Delta\spr
m.+\Gamma\spr m-1.$, equals
$$-2\int\limits_{W^{m}}\Delta_{(m-1)}(\Delta\spr m.(\Gamma\spr
m-1.)^2+(\Delta\spr m.)^2\Gamma\spr m-1.)=$$
$$-2\int\limits_{W^{m-1}}\Delta_{(m-1)}((\Gamma\spr m-1.)^2+
\Gamma\spr m-1.(-\sum\limits_{k=1}^{m-1}\omega\spr k.
+2\sum\limits_{k=2}^{m-1}\Delta\spr k.))$$ The third term III
clearly equals
$$\int\limits_{W^{m-1}}\Delta_{(m-1)}(\Gamma\spr m-1.)^2$$
The first term I equals
$$\int\limits_{W^m}\Delta_{(m)}w^*(\Gamma\spr m.)^2$$ so by the
projection formula it equals
 $$
 \int\limits_{X\sbr m._B}\Gamma_{(m)}(\Gamma\spr m.)^2=
 \int\limits_{\Gamma_{(m)}}(\Gamma\spr m.)^2.$$
 By Corollary 6, we have
 $$-\Gamma\spr m.|_{\Gamma_{(m)}}=e_m+\binom{m}{2}\omega$$
 where $e_m$ is the exceptional divisor of the natural map
 $\Gamma_{(m)}\to X$ (i.e. the restriction of the cycle map
 $\frak c_m$) and $\omega=\omega_{X/B}.$  By
 the discussion in \S 1.6, we have $e_m^2=-\beta_m.$,
 thus finally
 $$\int\limits_{W^m}\Delta_{(m)}w^*(\Gamma\spr m.)^2=
 -\sigma\beta_m+\binom{m}{2}^2\omega^2\tag 2.7$$
 where $\sigma$ as usual denotes the number of fibre nodes.
 This completes the integration of Case 1.\nl
 {\proclaim{Case 2}} $$B_1=c_1\spr i., B_2=c_2\spr j..$$
 This case is analogous to but much simpler than Case 1.
 We use the projection formula with respect to the ordered cycle map
 $p^m:W^m\to X^m_B$, which maps $\Delta_{(m)}$ with degree 1 to
 the small diagonal $X\subset X^m_B$,  and we
 get
 $$\int\limits_{W^m}\Delta_{(m)}c_1\spr i.c_2\spr j.=
 \int\limits_Xc_1c_2.\tag2.8$$
 {\proclaim{Case 3}} $$B_1=\Delta\spr i., B_2=c\spr j..$$
If $i,j<m$, we have as before
$$\int\limits_{W^m}\Delta_{(m)}\Delta\spr i.c\spr j.=
\int\limits_{W^{m-1}}\Delta_{(m-1)}\Delta\spr i.c\spr j.\tag 2.9$$
If $i<j=m,$ we have
$$\int\limits_{W^m}\Delta_{(m)}\Delta\spr i.c\spr m.=
\frac{1}{m-1}\sum\limits
_{k=1}^{m-1}\int\limits_{W^{m-1}}\Delta_{(m-1)}\Delta\spr i.c\spr
k.\tag 2.10$$ If $j<i=m$ we use (2.5) as before to obtain
$$\int\limits_{W^m}\Delta_{(m)}\Delta\spr m.c\spr j.=\frac{1}{m-1}
\int\limits_{W^{m-1}}\Delta_{(m-1)}(-\sum\limits_{k=1}^{m-1}\omega\spr
k. +2\sum\limits_{k=2}^{m-1}\Delta\spr k.)c\spr j.\tag 2.11$$
Finally if $i=j=m$ we may work as in Rule 3 above, representing
$c$ by a cycle disjoint from sing$(\pi)$, to conclude
$$\int\limits_{W^m}\Delta_{(m)}\Delta\spr m.c\spr m.=
-(m-1)\int\limits_X\omega.c\tag 2.12$$

 Putting together Rules 1-3 and Cases 1-3 above,
 we have completed  the integration of all top-degree
 monomials.\qed
 \remark{Remark 2.1} We note that the 'physical' small diagonal in
 $W^m$
 $$\Delta_{(m)}=w\inv(\Gamma_{(m)})=\Delta\spr2.\cap...\cap\Delta\spr m.$$
 will in general have 'extraneous', possibly improper-dimensional
  components $K$ resulting from
 the fact, noted in \cite{R1},
  that a punctual scheme belonging the component $C^m_i$ admits
 as subscheme an arbitrary punctual length-$(m-2j)$ belonging to
 $C^{m-2j}_i,$ provided $1\leq i\leq m-2j-1.$
This gives rise to a component of $\Delta_{(m)}$ birational
 to $\prod\limits_{j=1}^{[\frac{m-i-1}{2}]}C^{m-2j}_i$. Any such
 component $K$, however, automatically collapses under $w$,
 as it cannot dominate the irreducible $\Gamma_{(m)}$. Therefore
 $K$ does
 not contribute to $w_*(\Delta_{(m)})$
\endremark

\endcomment
\remark{Remark 2.5.1} Let $$u:E_1\to E_0$$ be a map of rank-2
vector-bundles on a scheme $X$, which drops rank by 1 along a
divisor $Z$, i.e. is locally of the form diag$(1,z)$, where $z$ is
an equation of $Z$. Then $u$ induces a rational map, known as an
'elementary modification'
$$\P(E_0)\dashrightarrow\P(E_1)$$ which is defined by a correspondence $$\matrix
&&Q&&\\&\alpha\swarrow&&\searrow\beta&\\\P(E_0)&&&&\P(E_1)\endmatrix$$
where $Q\subset \P(E_0)\times_X\P(E_1)$ is the 0-locus of the
natural map induced by $u$
$$p_2^*(M_{E_1})\to p_1^*(\O_{\P(E_0)}(1))$$ where $M_{E_1}$ is
the tautological subbundle (which in this case coincides with
$\O_{\P(E_1)}(-1)$ because $E_1$ has rank 2). Then
$$\beta^*(\O_{\P(E_1)}(1))=\alpha^*(\O_{\P(E_0)}(1))(-Z).\tag
2.5.1.1$$ Indeed (2.5.1.1) is obvious because by $Q$'s definition
there is a natural map induced by $u$,
$\beta^*(\O_{\P(E_1)}(1))\to\alpha^*(\O_{\P(E_0)}(1))$ and this has
divisor of zeros precisely $Z$.

\endremark
 \subheading{2.6 Example } With $X/B$ as above ($B$ a smooth curve),
 suppose $f:X\to\P^{2m-1}$ is a
 morphism. One, quite special, class of examples of this situation
 arises as what we call a {\it{generic rational pencil}}; that is,
 generally,
 the normalization of the family of rational curves of fixed degree $d$
  in $\P^r$ (so
 $r=2m-1$ here) that are incident to a generic collection $A_1,
 ...A_k$ of linear spaces, with
 $$(r+1)d+r-4=\sum(\text{codim}(A_i)-1);$$ see \cite{R3} and
 references therein, or \cite{RA} for an 'executive summary'.
 Then one expects a finite number $N_m$ of curves $f(X_b)$ to
 admit an m-secant $(m-2)$-plane, and this number can be evaluated as
 follows. Let $G=G(m-1,2m)$ be the Grassmannian of $(m-2)$-planes
  in $\P^{2m-1}$,
 with rank-$(m+1)$ tautological subbundle $S$, and let $L=f^*\O(1).$
 Then $$m!N_m=
 \int\limits_{W^m\times G}c_{m(m+1)}(S^*\boxtimes w^*\lambda_m(L))$$
 $$=\int\limits_{W^m\times G}c_{m+1}(S^*(L\spr 1.))
 c_{m+1}(S^*(L\spr 2.-\Delta\spr 2.))
 \cdots c_{m+1}(S^*(L\spr m.-\Delta\spr m.))$$
 $$=\int\limits_{W^m\times G}
 \prod\limits_{i=1}^m(\sum\limits_{j=0}^{m+1}\binom{m+1}{j}
 c_{m+1-j}(S^*)(L\spr
 i.-\Delta\spr i.)^j)$$
 $$=\sum\limits_{|(j.)|=m+1}\int\limits_Gc_{m+1-j_1,...,m+1-j_m}(S^*)
 \int\limits_{W^m}(L\spr
 1.)^{j_1}(L\spr 2.-\Delta \spr 2.)^{j_2}...(L\spr m.-\Delta\spr
 m.)^{j_m}$$
where $c_{u,v,w}=c_uc_vc_w$. Note that only terms with $j_m>0$
contribute. By the intersection calculus developed above, this
number can be computed in terms of the characters $L^2, \deg_\pi(L),
\omega^2, \sigma, \omega.L, \deg_\pi(\omega)=2g-2, g=$fibre genus;
in the generic rational pencil case, all these characters can be
computed by recursion on $d$.\par Suppose now that $m=3$, where the
only relevant $(j.)$ are
$$(2,1,1), (1,1,2), (1,2,1),(1,0,3),
(0,3,1),(0,2,2),(0,1,3),    (0,0,4).$$
 In each of these cases, it is easy to see that the $G$
 integral evaluates to 1. The $W$ integrals
 %which we will label respectively $I(j.)$,
 may be
evaluated by the calculus developed above. The relevant formulae are
$$\int\limits_{W^3}u\Delta\spr 3.=2\int\limits_{W^2}u, u\in T^2\tag
2.6.0$$
$$(\Delta\spr2.)^2=(\Gamma\scl 2.)^2=F^{(12:)}+q_{12}[\omega]\tag 2.6.1$$
(as usual we use $F^{(12:)}$ as short for
$F^{(12:\emptyset|\emptyset)}$)
$$\int_{W^2}L\spr i.(\Delta\spr 2.)^2=L.\omega=1/2\int\limits_{W^3}
L\spr i.(\Delta\spr 2.)^2\Delta\spr 3.,
i=1,2\tag2.6.2$$$$=1/2\int\limits_{W^3} L\spr 3.(\Delta\spr
2.)^2\Delta\spr 3., \tag 2.6.3$$
$$\int\limits_{W^2}L\spr i.L\spr j.\Delta\spr 2.=L^2=1/2\int\limits_{W^3}
L\spr i.L\spr j.\Delta\spr 2.\Delta\spr 3., (i,j)=(1,1), (1,2),
(2,2)\tag 2.6.4$$$$=1/2\int\limits_{W^3} L\spr i.L\spr 3.\Delta\spr
2.\Delta\spr 3., i=1,2,3;\tag 2.6.5$$
$$\int\limits_{W^2}L\spr 1.(L\spr
2.)^2=\deg_\pi(L)L^2=1/2\int\limits_{W^3}(L\spr 1.)L\spr 2.L\spr
3.\Delta\spr 3.=$$$$=\int\limits_{W^3}(L\spr 1.)^i(L\spr 2.)^j(L\spr
3.)^k\Delta\spr 3., (i,j,k)= (1,0,2), (0,1,2)\tag 2.6.6$$
$$\int\limits_{W^2}(\Delta\spr 2.)^3=-\sigma+\omega^2=
1/2\int\limits_{W^3}(\Delta\spr 2.)^3\Delta\spr 3.\tag 2.6.7$$
$$(\Delta\spr
3.)^2=2q_{123}[1]-q_{13}[\omega]-q_{23}[\omega]+F^{(13:)}+F^{(23:)}
\tag 2.6.8$$ where $F_s^{(i3:)}=\P(\O(-n_s))$ over $X'_s\coprod
X"_s$, with the 'correct' $\O(1)$, i=1,2;
$$L\spr 3..(\Delta\spr
3.)^2=2q_{123}[L]-q_{13}[\omega.L]-q_{23}[\omega.L]$$
%$$\int\limits_{W^3}(L\spr 3.)^2.(\Delta\spr
%3.)^2=2L^2$$
$$\matrix \int\limits_{W^3}L\spr 3.L\spr i..(\Delta\spr
3.)^2&=&2L^2-\deg_\pi(L)L.\omega,\ \  i=1,2\\ &=&2L^2,\ \
i=3\endmatrix\tag 2.6.9$$

$$\int\limits_{W^3}u(\Delta\spr
3.)^2=\int\limits_{W^2}u(2\Delta\spr 2.-\omega\spr 1.-\omega\spr
2.), u=L\spr 1.\Delta\spr 2.=L\spr 2.\Delta\spr 2., (\Delta\spr
2.)^2$$

(we can ignore $F$ terms because $u$ is perpendicular to them by
(2.3.5))
$$=\int\limits_{W^2}L\spr 1.(2q_{12}[-\omega]
-q_{12}[\omega]-q_{12}[\omega])$$
$$=-4L\omega,\ \text{ if}\  u=L\spr 1.\Delta\spr 2.=L\spr 2.\Delta\spr 2.\tag 2.6.10$$
$$=\int_{W^2}(\Delta\spr 2.)^2(2\Delta\spr 2.-\omega\spr
1.-\omega\spr 2.)=\int\limits_{W^2}2(\Delta\spr
2.)^3+2q_{12}[\omega^2]$$
$$=-2\sigma+4\omega^2,\ \text{if}\ u=(\Delta\spr 2.)^2;\tag 2.6.11$$
$$\matrix (\Delta\spr 3.)^3=&&2(\Gamma\scl 3.-\Gamma\scl
2.)q_{123}[1]\\&-2q_{123}[\omega]&+q_{13}[\omega^2]+q_{23}[\omega^2]+(\Gamma\scl
3.-\Gamma\scl 2.)(F^{(13:)}+F^{(23:)})\endmatrix\tag 2.6.12$$
$$\int\limits_{W^3}L\spr i.(\Delta\spr 3.)^3=2L.\omega, i=1,2\tag 2.6.13$$
$$\int\limits_{W^3}\Delta\spr 2.(\Delta\spr
3.)^3=-6\sigma+8\omega^2\tag 2.6.14$$
$$\int\limits_{W^3}(\Delta\spr
3.)^4=2(-3\sigma+4\omega^2)+2.2.\omega^2+\omega^2+\omega^2+2(-2\sigma+4\sigma)
$$$$=-2\sigma+14\omega^2\tag 2.6.15$$ where we have used the facts
$$(\Gamma\scl 3.)^2.\Gamma_{(3)}=\Gamma\scl
3.(F^{(123:)}+3q_{123}[\omega])=-6\sigma+9\omega^2, (\Gamma\scl
2.)^2.\Gamma_{(3)}=-\sigma+\omega^2,$$(both by \S1.6, as
$\beta_{3,1}=\beta_{3,2}=3, \beta_{2,1}=1$)$$ \Gamma\scl
3.\Gamma\scl 2.\Gamma_{(3)}=\frac{1}{2}(\Gamma\scl 3.)^2(\Gamma\scl
2.)^2=\frac{1}{2}\int\limits_{W^3}(\Gamma\scl
3.)^2\gamma^{3,2*}(F^{(12:)}+q_{12}[-\omega])$$$$=
\frac{1}{2}\int\limits_{F^{(12:)}}(\Gamma\scl
3.)^2+\frac{1}{2}\int\limits_{W^3}\Gamma\spr
3.(q_{12}[\omega^2]+2q_{123}[-\omega])$$$$=\half\int\limits_{F^{(12:)}}(\Gamma\scl
2.-2.\text{fibre})^2+
\half\int\limits_{W^3}2q_{123}[\omega^2]+2q_{123}[(-\omega).(-\omega)]
$$$$ =-2\sigma+3\omega^2$$
(for the last equality, note that $F^{(12:)}$ is a single point on
$W^2$ so  $(\Gamma\scl2.)^2=0$ on $F^{(12:)}$ and likewise on its
pullback on $W^3$);
$$(\Gamma\scl
3.)^2.F^{(i3:)}=-2\sigma, (\Gamma\scl 2.)^2.F^{(i3:)}=0, \Gamma\scl
3.\Gamma\scl 2..F^{(i3:)}=-2\sigma, i=1,2.$$ From all these, the
evaluation of $N_3$ is routine.

%\vfill\eject


\Refs \widestnumber\key{DSS}\ref\key A\by B. Ang\'eniol\book
Familles de Cycles Alg\'ebriques- Sch\'ema de Chow\finalinfo
Lecture Notes in Math. no. 896\publ Springer
\endref% \ref\key DS\by S. Donaldson, I. Smith\paper Lefschetz
%pencils and the canonical class for symplectic
%4-manifolds\finalinfo arXiv:math.SG/0012067\endref
\ref\key EG\by G. Ellingsrud, L. G\"ottsche\paper Hilbert schemes
of points and Heisenberg algebras\finalinfo ICTP lectures, 1999
(available at http://ictp.trieste.it) \endref
 \ref\key L\by
M. Lehn\paper Chern classes of tautological sheaves on Hilbert
schemes of points on surfaces\jour  Invent. math\vol 136\pages
157-207\yr 1999\endref \ref\key LS\by M. Lehn, C. Sorger\paper The
cup product of Hilbert schemes for K3 surfaces\jour Invent.
math\vol 152\pages 305-329\yr 2003\endref
 \ref\key N\by H. Nakajima\paper Heisenberg
algebra and Hilbert schemes of points on projective surfaces\jour
Ann. Math. \vol 145 \yr 1997\pages 379-388\endref
 \ref\key R\by Z. Ran \paper Geometry on
nodal curves (math.AG/0210209; update at\newline
http://math.ucr.edu/$\ \tilde{\ }$ ziv/papers/geonodal.pdf)
\endref\ref\key R1\bysame \paper A note on Hilbert schemes
of nodal curves (preprint available at\newline
http://math.ucr.edu/$\ \tilde{\ }$ ziv/papers/hilb.pdf or at
arXiv.org/math.AG/0410037)
\endref\ref\key R2\bysame\paper
 Cycle map on Hilbert Schemes of Nodal Curves
 (arXiv.org/math.AG/0410036)\endref
 \ref\key R3\bysame\paper The degree of the divisor of jumping
rational curves\jour Quart. J. Math.\yr 2001
 \pages 1-18\endref\ref\key RA\bysame\paper Rational
curves in projective spaces (notes available at
\newline http://math.ucr.edu/$\ \tilde{
}$ziv/papers/ratcurv.pdf)\endref

\ref\key Se\by E. Sernesi\book Topics on families of projective
varieties\finalinfo (Queens papers in pure and applied Math. vol.
73)\publ Queens Univ.\yr 1986\endref

%\ref\key Sm\by I. Smith\paper Serre-Taubes duality for
%pseudo-holomorphic curves \finalinfo arXiv:math.SG/0106220\endref


\endRefs
\enddocument